\documentclass[leqno,a4paper,titlepage,twoside,11pt]{article}

\usepackage[T1]{fontenc}%%<<ouvrez,>>fermez les guillemets%%

\usepackage[francais]{babel}

\usepackage{bbm,pifont,eufrak,eucal,mathrsfs,latexsym,graphicx,epsfig,makeidx,color,rotating,multirow}

\usepackage[all]{xy}
\xyoption{2cell}

\newcounter{infra}[page]

\newenvironment{dem}[1][]{%
{\bf D\'emonstration #1 : }}{%
\hspace*{\fill}\nolinebreak[1]\hspace*{\fill}\underline{\bf Q.e.d.}\\}

\newenvironment{dem*}[1][]{% 
{\bf D\'emonstration #1 : }}{%
 }

\newenvironment{eq*}{\begin{eqnarray*}}{\end{eqnarray*}}

\newtheorem{defi}{D\'efinition}

\newtheorem{thm}{Th\'eor\`eme}[section] \newtheorem{lem}[thm]{Lemme}
\newtheorem{pro}[thm]{Proposition} \newtheorem{cor}{Corollaire}[thm]

\newcommand{\croi}{\times}

\newcommand{\A}{\ens A} 
\newcommand{\adh}{\overline}

 \newcommand{\boufl}{\xymatrix{ \!\!
\ar@(ur,dr) }}%%%%%boucle fl\'ech\'ee%%%%%

\newcommand{\CCC}{{\adh{{\cal C}}}}
 
\newcommand{\cad}{\mbox{\it c-\`a-d }}
\newcommand{\call}{\mathscr}

\newcommand{\cd}{\rangle}
\newcommand{\cf}{{\it cf. }}

\newcommand{\cg}{\langle}

\newcommand{\ch}{\vee} 
 
\newcommand{\chap}{\widehat}
\newcommand{\cl}{\mathrm{cl}}
 \newcommand{\codim}{\mathrm{codim}}

 \newcommand{\con}{\supseteq}

 \newcommand{\cten}{\Box
\hspace*{-1.77ex}\raisebox{.19ex}{$\croi$}}

 \newcommand{\demi}{\frac{1}{2}}

\newcommand{\dij}{\sqcup} 
\newcommand{\Dij}{\bigsqcup}

\newcommand{\donne}{\mapsto}

 \newcommand{\ens}{\mathbbm}

\newcommand{\equi}{\Leftrightarrow} \newcommand{\et}{\mathrm{\; et \;}}

\newcommand{\exist}{\exists\:}

\newcommand{\flfl}{\raisebox{-1ex}{$\check{}$}$\hskip -.2ex\frown$\hskip -.2ex\raisebox{-1ex}{$\check{}$}}
\newcommand{\Gr}{\mathrm Gr}

  \newcommand{\GL}{\mathrm{GL}} \newcommand{\goth}{\mathfrak}

\newcommand{\Id}{\mathrm{Id}}
 \newcommand{\ie}{\emph{i.e. }}
\newcommand{\im}{\mathrm{Im\:}} \newcommand{\impliq}{\Rightarrow}

 \newcommand{\infi}{\infty}

\renewcommand{\int}{\mathrm{int\:}}

 \newcommand{\inv}{^{-1}}
 \newcommand{\iso}{\simeq}

\newcommand{\kk}{{\mathbf{k}}}
\newcommand{\K}{\mathbbm{K}}

\newcommand{\loccit}{\emph{loc. cit.}}

\newcommand{\ma}{\displaystyle} %%%%%%matrices d'ordre 2, 3%%%%%%%

\newcommand{\moins}{\:\setminus\:} %mathenrelief%

 \newcommand{\Ou}{\mbox{\Ou}}

\newcommand{\p}{\:\:.}  
\newcommand{\PGL}{\mathrm{PGL}}
\newcommand{\PP}{\mathbbm{P}}

\newcommand{\Pic}{\mathrm{Pic}\;} 
\newcommand{\pic}{\mathrm{pic}\;}
\newcommand{\plus}{\oplus} \newcommand{\Plus}{\bigoplus}

\newcommand{\PSO}{\mathrm{PSO}}

\newcommand{\qq}{\forall\:}

\newcommand{\res}[1]{{\left | {}_{#1} \right.}}

\newcommand{\rang}{\mathrm{rang}}

\newcommand{\rg}{\mathrm{rg}\,}

\newcommand{\si}{\mathrm{\;si\;}} 
\newcommand{\sinon}{\mbox{ sinon}} 
\newcommand{\SL}{\mathrm{SL}}
\newcommand{\SO}{\mathrm{SO}}

\newcommand{\Sp}{\mathrm{Sp}}
  
\newcommand{\sub}{\subseteq}
\newcommand{\substr}{\raisebox{-1.5ex}{$\:\stackrel{\displaystyle
\subset}{\scriptstyle\not=}\:$}}

\newcommand{\Sym}{\mathrm{Sym}}
\newcommand{\Tens}{\bigotimes} \newcommand{\tens}{\otimes} 

\newcommand{\tenso}[1]{\raisebox{-1.5ex}{$\ma \stackrel{\displaystyle
\tens}{\scriptstyle #1}$}} 

\newcommand{\tilda}{\widetilde} 

\newcommand{\tq}{\: : \:}

\newcommand{\Uni}{\bigcup}

\newcommand{\uni}{\cup}

\newcommand{\Vect}{\bigwedge\nolimits}

\newcommand{\vide}{\emptyset}

\newcommand{\Z}{\ens Z} 

\newcommand{\zzz}{{\mathbf{z}}}

\newcommand{\appli}[5][]{\begin{array}{ccc} #2 &
\stackrel{#1}{\longrightarrow} & #3 \\ #4 & \longmapsto & #5
\end{array}}

\makeindex

\begin{document}

\title{Sur la cohomologie \`a support des fibr\'es en droites sur les vari\'et\'es sym\'etriques complètes}
 \author{Alexis TCHOUDJEM\\
  Institut Camille Jordan\\
Universit\'e Claude Bernard Lyon I\\
Boulevard du Onze Novembre 1918\\
69622 Villeurbanne\\
FRANCE
} \date{Villeurbanne, le \today}

\maketitle

{\bf R\'esum\'e : } \'Etant donn\'e un groupe alg\'ebrique lin\'eaire semi-simple $G$, on s'int\'eresse aux compactifications magnifiques des $G-$espaces homogènes sym\'etriques. Si $X$ est une telle compactification, si $L$ est un fibr\'e en droites $G-$lin\'earis\'e sur $X$ et si $C$ est une cellule de Bialynicki-Birula de $X$ de codimension $c$, alors l'algèbre  de Lie $\goth g$ de $G$ opère naturellement sur le groupe de cohomologie \`a support $H^c_C(L)$. On donne ici une condition n\'ecessaire, portant sur la cellule $C$, pour que ce $\goth g-$module possède un sous-quotient simple de dimension finie. On en d\'eduit une formule pour la caract\'eristique d'Euler-Poincar\'e de $L$ sur $X$ et une estimation (exacte pour certains cas dont celui de la vari\'et\'e des coniques complètes) des groupes de cohomologie sup\'erieure $H^d(X,L)$, $d \ge 0$.

\vskip 1cm

\newpage

\tableofcontents

\vskip .5cm

Soit $ \kk$ un corps alg\'ebriquement clos de caract\'eristique nulle.

\addcontentsline{toc}{section}{Introduction}

\section*{Introduction}

Consid\'erons un groupe alg\'ebrique lin\'eaire semi-simple $G$ sur $\kk$ et $\goth g$ son algèbre de Lie. Si $X$ est une $G-$vari\'et\'e projective lisse et $L$ un fibr\'e en droites $G-$lin\'earis\'e sur $X$, alors tous les groupes de cohomologie $H^d(X,L)$, $d \ge 0$, sont des $G-$modules de dimension finie.

Lorsque $X$ possède une d\'ecomposition cellulaire :
$$X = \Dij_i C_i$$
où les cellules $C_i$ sont des sous-vari\'et\'es de $X$, de codimension $c_i$, isomorphes \`a des espaces affines, on dispose d'une suite spectrale, la {\it suite spectrale de Grothendieck-Cousin} (\cf \cite{Kempf}), dont les premiers termes sont les groupes de cohomologie \`a support $H^{c_i}_{C_i}(L)$ et qui converge vers les groupes de cohomologie globaux $H^d(X,L)$.

Remarquons que les $H^{c_i}_{C_i}(L)$ sont seulement des $\goth g-$modules alors que les $H^d(X,L)$ sont des $G-$modules.

Dans certains cas, la partie finie de la suite spectrale de Grothendieck-Cousin d\'eg\'enère, \ie les termes de la suite spectrale ont des multiplicit\'es constantes selon chaque $\goth g-$module simple de dimension finie. Cela permet alors de d\'eterminer compl\`etement les $H^d(X,L)$. Cela se produit si $X$ est une vari\'et\'e de drapeaux (d'où le c\'elèbre th\'eorème de Borel-Weil-Bott) et plus g\'en\'eralement si $X$ est une vari\'et\'e magnifique de rang minimal (ce sont des $G-$vari\'et\'es projectives avec une seule $G-$orbite ferm\'ee et qui v\'erifient notamment la prpopri\'et\'e d'avoir une d\'eciomposition cellulaire dont toutes les cellules rencontrent cette orbite ferm\'ee), \cf \cite[th. 3.1]{toi2}.
On s'int\'eresse ici aux $G-$vari\'et\'es sym\'etriques complètes d\'efinies par De Concini et Procesi dans \cite{DCP}. Ces vari\'et\'es ont aussi une seule orbite ferm\'ee.
Lorsqu'une cellule $C_i$ rencontre cette orbite projective, le $\goth g-$module $H^{c_i}_{C_i}(L)$ possède une suite de composition finie dont les quotients successifs sont des $\goth g-$modules bien \'etudi\'es : les modules de Verma tordus (\cf \cite[th. 4.1]{toi}).

Mais quand une cellule $C_i$ ne rencontre pas la $G-$orbite projective de $X$, on n'a pas de description aussi commode du $\goth g-$module $H^{c_i}_{C_i}(L)$. N\'eanmoins, nous allons montrer (avec une condition suppl\'ementaire sur le fibr\'e en droites $L$) que, dans ce cas, le $\goth g-$module $H^{c_i}_{C_i}(L)$ n'a pas de $G-$modules parmi ses sous-quotients simples (c'est le th\'eorème \ref{th:orbf}).
Ce r\'esultat permet de simplifier la suite spectrale de Grothndieck-Cousin. Comme applications, on donnera une formule pour la caract\'eristique d'Euler-Poincar\'e global (\cf le th\'eorème \ref{thm:euler}) :
$$\sum_i (-1)^iH^i(X,L)$$
(somme altern\'ee dans l'anneau des repr\'esentations de $G$), une majoration des multiplicit\'es des $G-$modules simples dans les $H^d(X,L)$ (\cf le th\'eorème \ref{thm:ineg}) avec, pour corollaire un r\'esultat d'annulation en degr\'e $d=1$ (\cf le corollaire \ref{cor:annul1}).

Enfin on d\'eterminera complètement les $H^d(X,L)$ pour deux exemples de vari\'et\'es sym\'etriques complètes de rang $2$ dont la vari\'et\'e des coniques complètes (\cf les th\'eorèmes \ref{th:chq} et \ref{thm:ccc}).

\printindex

\section{\' Enonc\'e du r\'esultat principal}

\subsection{Vari\'et\'es sym\'etriques complètes}

Soit $G$ un groupe alg\'ebrique lin\'eaire sur $\kk$. On suppose que $G$ est semi-simple adjoint et connexe sur $\kk$. Soit $\theta : G \to G$ un automorphisme d'ordre $2$. Notons $H:= G^\theta$ le sous-groupe des points fixes de $\theta$.

L'espace homogène $G/H$ est une {\it vari\'et\'e sym\'etrique} affine. 

D'après \cite[th. 3.1]{DCP}, il existe une unique $G-$vari\'et\'e projective $X$ telle que :

\begin{description}
\item[i)] $X$ contient $G/H$ comme $G-$orbite ouverte ;
\item[ii)] $X$ est lisse ;
\item[iii)] le diviseur $X \moins G/H = D_1 \cup ... \cup D_r$ est un diviseur \`a croisements normaux (les composantes irr\'eductibles $D_i$ de $X \moins G/H$ sont appel\'ees les diviseurs limitrophes de $X$) ;
\item[iv)] pour tous $x,x' \in X$, $G.x=G.x' \equi \left\{ i \tq x \in D_i \right \} = \left\{ i \tq x' \in D_i \right \}$ ;
\item[v)] il n'y a qu'une seule $G-$orbite ferm\'ee dans $X$ : $F=D_1 \cap ... \cap D_r$. 

\end{description}
L'entier $r$ est le {\it rang} de $X$.

On dit que la vari\'et\'e $X$ est la {\it compactification magnifique} de $G/H$. On appelle {\it vari\'et\'es sym\'etriques complètes} de telles vari\'et\'es.

{\bf Remarque :} Soit $\goth h$ l'algèbre de Lie de $H$. L'automorphisme $\theta$ induit un automorphisme d'algèbres de Lie $\goth g \to \goth g$ encore not\'e $\theta$. On a :
$$\goth h = \goth g^\theta = \left\{ x \in \goth g \tq \theta(x) = x\right\} \p$$
Soit $h:= \dim \goth h$. La vari\'et\'e $X$ peut \^etre d\'efinie comme l'adh\'erence de l'orbite $G.\goth h$ dans la grassmannienne $\call{G}:=\mathrm{Gr}_h (\goth g)$ des sous-espaces de $\goth g$ de dimension $h$ (\cf \cite[\S 6]{DCP}).

Fixons pour la suite $G$, $H$ comme ci-dessus et $X$ la compactification magnifique de l'espace sym\'etrique $G/H$. On notera $\chap{G}$ le rev\^etement universel de $G$ et on consid\'erera $X$ comme une $\chap{G}-$vari\'et\'e.

{\bf Exemple :} Soit $G:= \PGL_3$ et soit $\theta$ l'involution :
$$\theta : G \to G \; , \; [g] \donne [(g^t)\inv ]$$
(on note $[g] \in \PGL_3$ la classe modulo $\kk^*$ d'une matrice  inversible $g$). On a dans ce cas $H=G^\theta=\SO_3$. 
Notons $S_3$ l'espace des matrices sym\'etriques $3 \croi 3$ \`a coefficients dans $\kk$.
Soit ${\cal C}$ la vari\'et\'e alg\'ebrique affine des coniques non d\'eg\'en\'er\'ees de $\PP^2$ (vues comme des classes modulo $\kk^*$ de formes quadratiques non d\'eg\'en\'er\'ees sur $\kk^3$, vues elles-m\^emes comme des matrices sym\'etriques $3 \croi 3$ non singulières \`a coefficients dans $\kk$) :
$${\cal C} := \{[q] \in \PP(S_3) \tq \det q \not=0\} \p$$

Comme le groupe $G$ agit transitivement sur ${\cal C}$, on peut identifiier l'espace homogène sym\'etrique $G/H$ \`a ${\cal C}$ via l'isomorphisme :
$$G/H \to {\cal C} \; , \; [g]H \donne [(g^t)\inv g] \p$$

On peut dans ce cas d\'ecrire la compactification magnifique $\CCC$ de ${\cal C}$ comme une sous-vari\'et\'e ferm\'ee de $\PP^5 \croi \PP^5$ par :
$$\CCC:= \{([q],[q']) \in \PP(S_3) \croi \PP(S_3) \tq qq' \in \kk \Id\}$$
où $qq'$ est le produit usuel des matrices $q,q'$.

En effet, $\CCC$ est lisse, contient ${\cal C}$ via l'inclusion :
$$[q] \in {\cal C} \donne ([q],[q\inv]) \in \CCC$$
et le groupe $\chap{G}=\SL_3$ agit sur $\CCC$ par :
$$\qq g \in \chap{G}, \qq  ([q],[q']) \in \CCC, \; g. ([q],[q']) = ([(g^t)\inv q g\inv] , [gq'g^t]) \p$$
 Pour cette action, les diviseurs limitrophes de $\CCC$ sont :
$$D_1:= \{([q],[q']) \in \CCC \tq \rg q =1\}\et D_2 :=\{([q],[q']) \in \CCC \tq \rg q' =1\}$$ et les axiomes de d\'efinition d'une vari\'et\'e $G-$magnifique sont bien v\'erifi\'es. Notons que $\CCC$ est de rang $2$ ; c'est la {\it vari\'et\'e des coniques complètes.}

{\bf Remarque :} On peut aussi d\'efinir $\CCC$ comme l'\'eclat\'e de $\PP^5$ le long de la surface de Veronese $\nu(\PP^2)$ où $\nu : \PP^2 \to \PP^5$, $[x:y:z]\donne[x^2:y^2:z^2:xy:xz:yz]$.

\subsection{Faisceaux inversibles sp\'eciaux}

Le r\'esultat principal de cet article concerne certains faisceaux inversibles sur $X$ :

\begin{defi}
On dit qu'un faisceau inversible $\call{L}$ sur $X$ est {\it sp\'ecial} s'il existe un entier $n >0$ tel que :
$$\call{L}^{\tens n} \iso {\cal O}_X(n_1D_1+...+n_r D_r)$$
pour certains entiers $n_i$.

\end{defi} 
{\bf Remarque :} si la vari\'et\'e $X$ n'est pas exceptionnelle \ie si le rang du groupe de Picard de $X$ est $r$, alors tous les faisceaux inversibles sur $X$ sont sp\'eciaux. C'est par exemple le cas pour la compactification magnifique de $PGL_n/PSO_n$ (pour tout $n \ge 2$).

\subsection{Cellules et cohomologie \`a support}

Soit $\nu$ un sous-groupe \`a un paramètre de $G$ tel que l'ensemble des points fixes :
$$X^\nu := \{x \in X \tq \qq s \in \kk^*, \nu(s).x = x\}$$
est fini. De  tels sous-groupes \`a un paramètre existent toujours et sont appel\'es {\it $X-$r\'egulier}.  

Pour tout point fixe $x \in X^\nu$, on pose :
$$X^+_x \mbox{ ou } X^+(x):=\{y \in X \tq \lim_{s \to 0 \atop s \in \kk^*} \nu(s).y = x \} \p$$

D'après \cite{BiBi}, chaque $X^+(x)$ est une sous-vari\'et\'e localement ferm\'ee de $X$, isomorphe \`a un espace affine. On obtient ainsi une d\'ecomposition : $X = \dij_{x \in X^\nu } X^+(x)$. Les $X^+(x)$ sont les {\it cellules de Bialynicki-Birula}  de $X$ et on dit que $x$ est le centre de la cellule $X^+(x)$.

Tous les faisceaux inversibles $\call{L}$ sur $X$ sont $\chap{G}-$lin\'earis\'es (\cf \cite{St}). Donc pour tout point $x \in X^\nu$ et tout entier $n$, les groupes de cohomologie \`a support $H^n_{X^+(x)}(\call{L})$ sont des $\goth g-$modules (\cf \cite[lem. 11.1]{Kempf}).

De plus, $H^n_{X^+(x)}(\call{L}) = 0$ sauf si la cellule $X^+(x)$ est de codimension $n$ dans $X$.

{\bf Exemple :} Soient $G:= \PGL_2(\kk)$ et $X$ la compactification magnifique de l'espace homog\`ene $\PGL_2/\PSO_2$. On a : $\chap{G}=\SL_2(\kk)$ et  la variété $X$ est l'espace projectif $\PP(V)=\PP^2$ où $V$ est le $\SL_2(\kk)-$module simple $\kk[T_0,T_1]_2$ des polyn\^omes homog\`enes de degré $2$. L'unique $G-$orbite fermée est formée des classes de polyn\^omes $[P]$ de discriminant nul. On a :
$$\chap{T} = \left\{{\left(\begin{array}{cc}
s & 0 \\
0 & s\inv
\end{array}\right) \tq s \in \kk^*}\right\} \et {\cal X} = \Z \omega \mbox{ où } \omega : \left(\begin{array}{cc}
s & 0 \\
0 & s\inv
\end{array}\right) \donne s \p$$

Si on prend pour sous-groupe à un param\`etre $\nu : \kk^* \to G $, $s \donne \left[\begin{array}{cc}
s & 0 \\
0 & s\inv
\end{array}\right] $, alors on a :
$$X^\nu =\left\{ [T_0^2] , [T_0T_1],[T_1^2]\right\} \p$$

Seul le point fixe $x_0 := [T_0T_1]$ n'est pas dans la $G-$orbite fermée. La cellule correspondante est :

$$X^+(x_0)=\left\{ [bT_0T_1 + c T_1^2] \in X \tq b ,c \in \kk, b\not=0\right\} \iso \A^1$$

  Les faisceaux inversibles sur $X$ sont les ${\cal O}_{\PP^2}(n )$ avec  $n$ entier.

Comme le point $x_0$ est $N_{\chap{G}}(\chap{T})-$stable et comme tous les caractères du groupe $N_{\chap{G}}(\chap{T})$ sont triviaux, le tore $\chap{T}$ agit trivialement sur la fibre ${\cal O}_{\PP^2}(n )\res{x_0}$ (pour chaque entier $n$). On a donc un isomorphisme de $\chap{T}-$modules :

$$H^1_{X^+(x_0)}({\cal O}_{\PP^2}(n )) \iso H^1_{X^+(x_0)}({\cal O}_X)$$
pour chaque $n \in \Z$. 

Or, les poids de l'espace tangent $T_{x_0}X$ sont $-2\omega$, $0$ et $2 \omega$. On en déduit que tous les poids du $\chap{T}-$module $H^1_{X^+(x_0)}({\cal O}_{\PP^2}(n ))$ sont de la forme $2n\omega$ avec $n \in \Z_{<0}$. En particulier, aucun n'est dominant. En conséquence, le $\mathrm{sl}_2-$module $H^1_{X^+(x_0)}({\cal O}_{\PP^2}(n ))$ ne peut pas avoir de $\SL_2(\kk)-$modules parmi ses sous-quotients simples.

Le théor\`eme qui suit est une généralisation de ce fait aux variétés symétriques compl\`etes :

\begin{thm}
Soient $\call{L}$ un faisceau inversible sp\'ecial sur $X$ et $C$ une cellule de Bialynicki-Birula de $X$ de codimension $d$.

Si le $\goth g-$module $H^d_C(\call{L})$ admet un sous-quotient simple de dimension finie, alors le centre de la cellule $C$ est dans l'orbite ferm\'ee de $X$.
\end{thm}

{\bf Remarque :} si $X$ est de rang minimal, \ie si $r=\rang(G)-\rang(H)$ (par exemple si $X$ est la compactification magnifique de $G \croi G / G$), alors tous les points fixes $x \in X^\nu$ sont dans $F$. Donc ce th\'eorème n'apporte rien  dans ce cas. En revanche si $X$ n'est pas de rang minimal (par exemple si $X$ est la compactification magnifique de $PGL_n/PSO_n$), alors il existe des points fixes $x \in X^\nu \moins F$.

Avant de d\'emontrer ce th\'eorème, nous allons rappeler quelques r\'esultats concernant les vari\'et\'es sym\'etriques complètes.

\section{Sous-groupe de Borel, tore maximal, système de racines, etc. adapt\'es \`a une vari\'et\'e sym\'etrique}\label{sec:debut}

Nous suivons \cite[\S 1]{DCP}.

Fixons un tore $T_1$ de $G$, {\it anisotrope} (\ie $\theta(t) = t\inv$ pour tout $t$ dans $T_1$) et maximal. Soit $T$ un tore maximal de $G$ qui contient $T_1$. Le tore $T$ est forc\'ement $\theta-$stable. On notera $\chap{T}$ l'image r\'eciproque de $T$ dans $\chap{G}$, ${\cal X}$ le r\'eseau des caractères de $\chap{T}$, $\cg \cdot, \cdot \cd$ le crochet naturel entre caractères et sous-groupes \`a un paramètre de $\chap{T}$. 

Soit $W$ le groupe de Weyl associ\'e \`a $(\chap{G},\chap{T})$. On choisit une forme bilin\'eaire $(\cdot,\cdot)$ sym\'etrique, non d\'eg\'en\'er\'ee et $W-$invariante
sur ${\cal X}$.

Si on note $\goth t$ l'alg\`ebre de Lie de $T$, on a :
$$\goth t = \goth t_0 + \goth t_1$$
où ${\goth t}_0:=\ker (\theta-1) \cap \goth t \et \goth t_1 := \ker (\theta +1) \cap \goth t$. 

Soit $\Phi \sub \goth t^*$ le système de racines de $(\goth g, \goth t)$. On note encore $\theta$ l'automorphisme induit par $\theta$ sur $\goth t^*$.  Cet automorphisme $\theta$ pr\'eserve $\Phi$ et aussi $\cg \cdot, \cdot \cd$ et $(\cdot , \cdot)$.

Posons $\Phi_0:=\left\{ \alpha \in \Phi \tq \theta(\alpha)=\alpha \right\}$ et $\Phi_1 := \Phi \moins \Phi_0$. On peut choisir l'ensemble des racines positives $\Phi^+$ de $\Phi$ tel que :
$$\qq \alpha \in \Phi^+ \cap \Phi_1, \; \theta(\alpha) \in \Phi^- \p$$

On notera $\ma \rho := \frac{1}{2}\sum_{\alpha \in \Phi^+}\alpha$. Ce $\rho$ est aussi un caractère de $\chap{T}$. 

Soient $\Delta$ la base d\'efinie par $\Phi^+$, $\Delta_0:= \Delta \cap \Phi_0$ et $\Delta_1:= \Delta \cap \Phi_1$.

Il existe une bijection $\adh{\theta} : \Delta_1 \to \Delta_1$ telle que :
$$\qq \alpha\in \Delta_1, \theta(\alpha) = -\adh{\theta}(\alpha) - \sum_{\delta \in \Delta_0} m_{\delta,\alpha}\delta$$
pour certains entiers $m_{\delta,\alpha} \ge 0$. Les entiers $m_{\delta,\alpha}$ sont entièrement d\'etermin\'es par les \'equations :
$$\cg \theta(\alpha) , \delta^\ch \cd = \cg \alpha , \delta^\ch \cd \; (\qq \alpha \in \Delta_1, \qq \delta \in \Delta_0)  \p$$

Pour toute racine $\alpha \in \Phi_1$, on pose $\tilda{\alpha} := \alpha - \theta (\alpha)$ et on num\'erote les racines de $\Delta_1$ : $\alpha_1,...,\alpha_d$ avec $d \ge r$, de sorte que :
$$\tilda{\Delta_1} := \left\{ \tilda{\alpha} \tq \alpha \in \Delta_1 \right\} = \left\{\tilda{\alpha_1},...,\tilda{\alpha_r}\right\} \p$$

\vskip .5cm

Enfin, on note $\chap{B}$ (respectivement $\chap{B}^-$) le sous-groupe de Borel de $\chap{G}$ d\'efini par $\Phi^+$ (respectivement par $\Phi^-$) et $B$ (respectivement $B^-$) son image dans $G$. 

\section{Groupe de Picard}

Nous rappelons ici la description du groupe de Picard de la vari\'et\'e $X$ comme  sous-r\'eseau du r\'eseau des poids de $\chap{T}$.

Puisque $X$ n'a qu'une seule $\chap{G}-$orbite ferm\'ee, il existe un unique point fixe $\zzz \in X$ du sous-groupe de Borel $\chap{B}^-$. Soit $Q$ le groupe d'isotropie de ce point $\zzz$ dans $\chap{G}$. C'est un sous-groupe parabolique de $G$ qui contient $\chap{B}^-$ et $F =\chap{G}.\zzz \iso \chap{G}/Q$.

Soit $ \call{L}$ un faisceau inversible sur $X$. Puisque le groupe $\chap{G}$ est simplement connexe le faisceau $\call{L}$ admet une unique $\chap{G}-$lin\'earisation (\`a isomorphisme près). On peut donc d\'efinir sans ambigu\"{i}t\'e le caractère $p(\call{L}) : Q \to \kk^*$ avec lequel  $Q$ agit sur la fibre $\call{L}\res{\zzz}$.

Si on note ${\cal X}(Q)$ le r\'eseau des caractères de $Q$, on a :
\begin{pro}[pro 8.1 de \cite{DCP}]\label{pro:inj}
Le morphisme $$\appli{\Pic(X)}{{\cal X}(Q)}{\call{L}}{p(\call{L})}$$
est injectif.
\end{pro} 

Notons $\pic(X)$ l'image de $\Pic(X)$ dans le r\'eseau ${\cal X}(Q)$ de sorte que : $\Pic (X) \iso \pic (X)$. 

D\'esormais, pour tout $\lambda \in \pic(X)$, on notera $\call{L}_\lambda$ un faisceau inversible $\chap{G}-$lin\'earis\'e sur $X$ de poids $p(\call{L}_\lambda)=\lambda$ et $[\call{L1}_\lambda]$ sa classe d'isomorphisme.

Quitte \`a renum\'eroter les diviseurs limitrophes $D_i$, $1 \le i \le l$, on supposera que : ${\cal O}_X(D_i) \iso \call{L}_{\tilda{\alpha_i}}$.

{\bf Exemple :} Dans le cas où $X=\CCC$, on a $\chap{G} = \SL_3$ et on peut prendre pour tore maximal $\chap{T}$ le tore des matrices diagonales de $\SL_3$. On a $\zzz= (\left[\begin{array}{ccc}
1 & 0 &0\\
0 & 0 & 0\\
0&0&0
\end{array}\right] , \left[\begin{array}{ccc}
0&0&0\\
0&0&0\\
0&0&1
\end{array}\right])$ et $Q=\chap{B^-}$, le sous-groupe des matrices triangulaires inf\'erieures de $\SL_3$. Si on note $\alpha_1,\alpha_2$ les racines simples d\'efinies par $\chap{G}$, $\chap{B}$ et $\chap{T}$ et $\omega_1,\omega_2$, les poids fondamentaux correspondants, on a :
$$\pic (\CCC) = 2\Z\omega_1 \plus 2\Z\omega_2$$
et on peut d\'ecrire explicitement les faisceaux inversibles sur $\CCC$ par :
$$[\call{L}_\lambda] = \left[{\cal O}_{\PP(S_3)}(m) \cten {\cal O}_{\PP(S_3)}(n) \res{\CCC}\right]$$
pour tout $\lambda = 2m \omega_1 + 2n\omega_2 \in 2\Z\omega_1 \plus 2\Z\omega_2$. 

\begin{center}
*
\end{center}

Dans la suite, on s'int\'eressera plus particulièrement aux faisceaux inversibles {\it sp\'eciaux} :

\begin{defi}
Un caractère $\lambda$ de $\chap{T}$ est {\it sp\'ecial} si $\theta(\lambda) = -\lambda$. 
\end{defi}

{\bf Remarque :} Soient $\lambda \in {\cal X}$ et $\call{L}_\lambda$ le faisceau inversible associ\'e. Le poids $\lambda$ est sp\'ecial si et seulement si le faisceau $\call{L}_\lambda$ l'est.

Soient $\omega_\alpha$, $\alpha \in \Delta$, les poids fondamentaux du système de racines $\Phi$ (consid\'er\'es comme des caractères de $\chap{T}$).

Voici une description de ces faisceaux inversibles sp\'eciaux \`a isomorphime près :

\begin{pro}[\S 2, th. 1 de \cite{CHMA}]
Le r\'eseau $$\left\{\lambda \in \pic(X) \tq \theta(\lambda)=-\lambda \right\}$$ admet pour base les poids $\tilda{\omega_1}, ..., \tilda{\omega_r}$ où :
$$\tilda{\omega_i} = \left\{ \begin{array}{ll}
\omega_{\alpha_i} + \omega_{\adh{\theta}(\alpha_i)} & \si \adh{\theta}(\alpha_i) \not= \alpha_i
\\
\omega_{\alpha_i} & \si \adh{\theta}(\alpha_i) = \alpha_i \et \theta(\alpha_i) \not= -\alpha_i \\
2 \omega_{\alpha_i} & \si \theta(\alpha_i) = -\alpha_i
\end{array}\right.$$
pour tout $1 \le i \le r$.
\end{pro}

{\bf Remarque :} Si $X$ est la compactification magnifique de $\PGL_n/\PSO_n$, on a $r=n-1$ et $\tilda{\omega_i}=2\omega_{\alpha_i}$ pour tout $1 \le i \le n-1$.

\section{Vari\'et\'es stables}

Les {\it vari\'et\'es stables} sont les sous-vari\'et\'es irr\'eductibles et $G-$stables de $X$ ; ce sont les vari\'et\'es $X_I := \bigcap_{i \in I} D_i$, $I$ partie de $\{1,...,r\}$.

Pour les vari\'et\'es stables, nous allons fixer quelques notations.

Pour toute partie $I \sub \{1,...,l\}$, on choisit un sous-groupe \`a un paramètre $\gamma_I: \kk^* \to T_1$ tel que pour tout $1 \le i \le r$ :
$$\cg \tilda{\alpha_i} , \gamma_I \cd >0 \si i \in I \;,$$
$$ \cg \tilda{\alpha_i}, \gamma_I \cd >0 \si i \not\in I \p$$

Soit $x_0:=H/H \in G/H \sub X$.

Si on pose $x_I:=\lim_{t\to 0} \gamma_I(t) . x_0$, alors on a $X_I=\adh{G.x_I}$.

{\bf Remarque :} Soient $I^c:= \{1,...,r\} \moins I$ et $\tilda{\Phi_{I^c}}$ les \'el\'ements de $\tilda{\Phi}$ qui sont combinaisons lin\'eaires (\`a coefficients entiers) des $\tilda{\alpha_i}$, $i \in I^c$.

Avec ces notations, $x_I$ est la sous-algèbre de Lie 
$$\goth t_0 \plus \Plus_{\alpha \in \Phi_0} \goth g_{\alpha} \plus \Plus_{\alpha\in \Phi_1 \atop {\tilda{\alpha} \in \tilda{\Phi}_{I^c}}} \kk .(X_\alpha + \theta(X_\alpha)) \plus \Plus_{\alpha \in \Phi_1^+ \atop \tilda{\alpha} \not\in \tilda{\Phi}_{I^c}} (\goth g_{\theta(\alpha)} \plus \goth g_{-\alpha}) $$
vue comme point de $X \sub \call{G}$.

On aura besoin aussi des groupes suivants :
$$P_I := \left\{ g \in G \tq \lim_{t \to 0} \gamma_I(t)\inv g \gamma_I(t) \mbox{ existe dans $G$ }\right\}$$
$$L_I := \left\{ g \in G \tq \qq t \in \kk^* , \gamma_I(t)\inv g \gamma_I(t) = g \right\} \p$$
Le groupe $P_I$ est un sous-groupe parabolique de $G$ contenant $Q \con B^-$ et $L_I$ est son sous-groupe de Levi par rapport \`a $T$.

Soit $\adh{L_I}$ le quotient de $L_I$ par son centre. L'automorphisme $\theta$ de $G$ induit un automorphisme $\adh{\theta_I}$ de $\adh{L_I}$ et on pose $\adh{H_I}:=\adh{L_I}^{\adh{\theta_I}}$.

Tout cela \'etant pos\'e, si $X(\adh{\theta_I})$ est la compactification magnifique de l'espace sym\'etrique $\adh{L_I}/\adh{L_I}^{\adh{\theta_I}}$, on a $X_I = G \croi^{P_I} X(\adh{\theta_I})$ ce qui signifie :
\begin{pro}[\cf \cite{DCP} \S 5]\label{pro:proj}
Il existe un morphisme surjectif et $G-$\'equivariant :
$$X_I \stackrel{\pi_I}{\to} G/P_I \;\; x_I \donne P_I/P_I$$
tel que $\pi_I\inv(P_I/P_I) \iso X(\adh{\theta_I})$ comme $\adh{L_I}-$vari\'et\'es.
\end{pro}

{\bf Remarque :} Il r\'esulte de cette proposition que $x_I$ est fix\'e par le tore $T$ si et seulement si $\adh{L_I}^{\adh{\theta_I}}$ contient le tore $\adh{T}$, l'image de $T$ dans $\adh{L_I}$.

Pour ce qui suit, on pose :
$$W^{L_I}:= \{w \in W \tq \qq \alpha \in \Phi_{L_I} \cap \Phi^+, w( \alpha) \in \Phi^+\}$$
($\Phi_{L_I}$ est l'ensemble des racines de $L_I$).
 
\section{Poids des faisceaux inversibles en les points fixes du tore}
Soit $x \in X^T$. Soit $I \sub \{1,...,l\}$ tel que $\adh{G.x} = X_I$. Comme le point $\pi_I(x)$ est un point fixe de $T$ dans $G/P_I$, il existe $w \in W^{P_I}$ tel que $\pi_I(x) = wP_I/P_I$.

Soient $\lambda \in \pic(X)$ et $\call{L}_\lambda$ le faisceau inversible et $\chap{G}-$lin\'earis\'e sur $X$ correspondant ; on note $\lambda_x$ le caractère avec lequel $\chap{T}$ agit sur la fibre $\call{L}_\lambda\res{x}$.  En g\'en\'eral, le point $x$ n'est pas un point fixe du sous-groupe parabolique $wP_Iw\inv$ mais on a :

\begin{lem}\label{lem:pfibre}
Pour tout poids sp\'ecial $\lambda \in \pic(X)$, le caractère $\lambda_x$ se prolonge en un caractère de $wP_Iw\inv$.
\end{lem} 

\begin{dem}
Comme le faisceau $\call{L}_\lambda$ est $\chap{G}-$lin\'earis\'e, nous allons seulement traiter le cas où $w=1$ \ie $\pi_I(x) = P_I/P_I$.  

\begin{center}
*
\end{center}

Supposons pour commencer que $\lambda$ est un caractère de $P_I$ \ie :
$$\qq \alpha \in \Phi, \cg \alpha,\gamma_I \cd = 0 \impliq (\lambda,\alpha) =0 \p$$
Dans ce cas, soit $\call{M}_\lambda$ le faisceau inversible induit par $\lambda$ sur la vari\'et\'e de drapeaux $G/P_I$. On note toujours  $\pi_I : X_I \to G/P_I$ la projection de la proposition \ref{pro:proj}. Si on note $\kk_\lambda$ la droite $\kk$ munie de l'action du tore $\chap{T}$ via le caractère $\lambda$, alors on a les isomorphismes de $\chap{T}-$modules :
$$\call{L}_\lambda\res{\zzz} \iso \kk_\lambda \iso \call{M}_\lambda\res{P_I/P_I} = \call{M}_\lambda\res{\pi_I(\zzz)} \iso \pi_I^*\call{M}_\lambda \res{\zzz} \p$$
On en d\'eduit d'après la proposition \ref{pro:inj} que $\call{L}_\lambda \iso  \pi_I^*\call{M}_\lambda $. En particulier, $\call{L}_\lambda \res{x} \iso \call{M}_\lambda\res{P_I/P_I}$ et le caractère $\lambda_x$ se prolonge en un caractère de $P_I$.
\begin{center}
*
\end{center}

Maintenant, ne supposons plus que $\lambda$ est un  caractère de $P_I$ :

Comme les $\tilda{\alpha_i}$, $1 \le i \le r$ sont $\Z-$lin\'eairement ind\'ependants, il existe des entiers $k_i$ , $i \in I^c$ et un entier $k >0$ tels que :
\begin{equation}\label{eq:gram}
\qq j \in I^c ,\; k(\lambda,\tilda{\alpha_j}) - \sum_{i \in I^c} k_i (\tilda{\alpha_i}, \tilda{\alpha_j}) =0 \p 
\end{equation} 

Posons $\mu := k\lambda - \sum_{i \in I^c} k_i \tilda{\alpha_i}$. C'est un poids sp\'ecial de $\pic(X)$.

De plus, pour toute racine  $\alpha \in \Phi$, on a :
$$\cg \alpha , \gamma_I \cd = 0 \impliq (\mu , \alpha)=0 \p$$
En effet, par exemple si $\alpha$ est une racine positive, $\alpha = \alpha_0 + \sum_{i=1}^d n_i \alpha_i$ pour un $\alpha_0 \in \Phi_0^+$ et certains entiers $n_i \ge 0$.

Donc $\tilda{\alpha} = \sum_{i=1}^d n_i \tilda{\alpha_i} = \sum_{i=1}^r n'_i \tilda{\alpha_i}$ pour certains entiers $n'_i$ tels que $n'_i \ge n_i$ pour tout $i$.

Or $\theta(\gamma_I) = - \gamma_I$ donc :
$$\cg \alpha , \gamma_I \cd = 0 \equi \cg \tilda{\alpha} , \gamma_I \cd = 0$$
$$\equi \sum_{i=1}^r n'_i  \cg \tilda{\alpha_i} , \gamma_I \cd = 0$$
\begin{equation}\label{eq:npiez}
\equi \qq i \in I, n'_i =0
\end{equation} 
car si $i \in I, \cg \tilda{\alpha_i} , \gamma_I \cd >0$.

Mais puisque $\theta(\mu) = - \mu$, on a pour chaque racine $\alpha$ :
$$(\mu,\alpha) = (\mu , \alpha_0) + \sum_{i=1}^d n_i(\mu , \alpha_i)$$
$$= \frac{1}{2} . \sum_{i=1}^d n_i (\mu , \tilda{\alpha_i}) $$
$$ = \frac{1}{2} . \sum_{i=1}^r n'_i (\mu , \tilda{\alpha_i})$$
$$ =  \frac{1}{2} . \sum_{i\in I} n'_i (\mu , \tilda{\alpha_i})$$
d'après (\ref{eq:gram}).

Par cons\'equent, on trouve gr\^ace \`a  (\ref{eq:npiez}) :
$$(\alpha,\gamma_I) = 0 \impliq (\mu,\alpha) = 0$$
autrement dit, $\mu$ est un caractère du sous-groupe parabolique $P_I$.

Rappelons que $\mu = k \lambda - \sum_{i \in I^c}k_i \tilda{\alpha_i}$. Ainsi :
$$\call{L}_\mu = \call{L}_{k\lambda} \tens \Tens_{i \in I^c} {\cal O}_X(-D_i)^{\tens k_i} \p$$

Or le point $x$ est dans la $G-$orbite ouverte de $X_I = \bigcap_{i \in I} D_i$ donc si $i \in I^c, x \not \in D_i$ et ${\cal O}_X(D_i) \res{x} = {\cal O}_X\res{x}$.
Par cons\'equent :
$$ \call{L}_\mu \res{x_I} \iso \call{L}_{k\lambda} \res{x} \p$$

Mais alors, le caractère $k \lambda_x$ avec lequel $\chap{T}$ agit sur $\call{L}_{k\lambda} \res{x}$ est un caractère de $P_I$ donc $\lambda_x$ aussi.
\end{dem}

\section{Espace tangent en un point fixe du tore}

Soit $x \in X^T$. Comme dans la section pr\'ec\'edente, on fixe $I \sub \{1,..,r\}$ tel que $X_I= \adh{G.x}$. On note encore $P_I$ le sous-groupe parabolique correspondant, $\pi_I : X_I \to G/P_I$ la projection associ\'ee et $L_I$ le sous-groupe de Levi de $P_I$ contenant $T$.  Soit $w \in W^{P_I}$ tel que $\pi_I(x) = wP_I/P_I$.

Notons $-\delta_1,...,-\delta_u$ les racines de $R^u(P_I)$, le radical unipotent de $P_I$.

{\bf Remarque :} les $\delta_j$ sont les racines positives de $\Phi$ telles que $\cg \delta_j , \gamma_I \cd \not= 0$.

L'espace tangent $T_x X_I$ est un $T-$module et on a :

\begin{lem}\label{lem:tang}
Les poids de $T_{x}X_I$ sont :
$$- w(\delta_1),...,-w(\delta_u), \pm w(\beta_1),...,\pm w(\beta_v) $$
pour certaines racines positives $\beta_1,...,\beta_v$ de $L_I$.
De plus, $2v=\dim X_I - \dim G/P_I$.
\end{lem}

{\bf Remarque :} Le groupe d'isotropie $G_x$ est contenu dans $wP_Iw\inv$ et $2v$ est aussi la codimension de $G_x$ dans $wP_Iw\inv$.

\begin{dem}
Il suffit de traiter le cas où $w=1$. Dans ce cas, $x \in \pi_I\inv(P_I/P_I) \iso \adh{L_I}/\adh{L_I}^{\adh{\theta_I}}$. En particulier le stabilisateur de $x$ dans $\adh{L_I}$ est un groupe r\'eductif car c'est un conjugu\'e de $\adh{L_I}^{\adh{\theta_I}}$. Notons $K$ ce stabilisateur.
D'après la proposition \ref{pro:proj}, on a un isomorphisme de $T-$modules :
$$T_{x} X_I \iso T_{P_I/P_I}G/P_I \plus T_{x} \adh{L_I}.x $$
$$\iso T_{P_I/P_I}G/P_I \plus T_{1} \adh{L_I}/T_{1}K  \p$$

D'une part, les poids de $T_{P_I/P_I}G/P_I$ sont les $-\delta_1,...,-\delta_u$. D'autre part, comme le groupe $K$ est r\'eductif, si $\beta$  est un poids de $T_{1} \adh{L_I}/T_{1} K$, $-\beta$ aussi. Donc les poids de $T_{1} \adh{L_I}/T_{1}K$ sont les $\pm \beta_1,..., \pm \beta_v$ où les $\beta_j$ sont les racines positives de $L_I$ qui ne sont pas racines de $K$.

 Enfin, on a aussi :
$$2v = \dim X_I -u = \dim X_I - \dim G/P_I$$
\end{dem}

\section{Th\'eorème principal}

\begin{thm}\label{th:orbf}
Soient $\call{L}$ un faisceau inversible sp\'ecial sur $X$ et $C$ une cellule de Bialynicki-Birula de $X$ de codimension $d$.

Si le $\goth g-$module $H^d_C(\call{L})$ admet un sous-quotient simple de dimension finie, alors le centre de la cellule, $x \in C^T$ est dans l'orbite ferm\'ee de $X$.
\end{thm}

Pour d\'emontrer ce th\'eorème on va calculer des caractères de $\chap{T}-$modules.

Rappelons d'abord les notations de \cite[\S 7.5.1]{Dix} et \cite{Kempf}:

Soit $M$ est un $\chap{T}-$module tel que pour tout caractère $\nu : \chap{T} \to \kk^*$, le sous-espace $M_\nu:=\{m \in M \tq \qq t \in \chap{T}, t.m = \nu(t)m\}$ est de dimension finie. Le caractère de $M$ est la fonction :
$$[M] : {\cal X} \to \Z \;,\; \nu \donne \dim M_\nu \p$$

 Soit une fonction $f : {\cal X} \to \Z$. On la notera   $\sum_{\nu \in {\cal X}} f(\nu) e^\nu$. Son {\it support} est l'ensemble $\{\nu \in {\cal X} \tq f(\nu) \not=0 \}$. 
Soit $\Z_{\ge 0} \Phi^+$ le c\^one de ${\cal X}$ engendr\'e $\Phi^+$. On notera $\Z\cg{\cal X}\cd$ l'anneau des fonctions $f : {\cal X} \to \Z$ dont le support est contenu dans une r\'eunion finie d'ensembles de la forme $\nu_i - \Z_{\ge 0}\Phi^+$, $\nu_i \in {\cal X}$.

{\bf D\'emonstration du th\'eorème :} 

On peut supposer que la cellule $C$ est de la forme :
$C =X^+(x)$ pour un certain point fixe $x \in X^T$ et un certain sous-groupe \`a un paramètre $\zeta$ dominant et r\'egulier (\ie $\cg \alpha ,\zeta \cd >0$ pour toute racine positive $\alpha$) tel que $X^{\zeta(\kk^*)} = X^T$).

Dans toute la suite on notera pour tout caractère $\nu \in {\cal X}$ :
$$|\nu| := \nu \si \cg \nu ,\zeta \cd \ge 0 \;, \; -\nu \si   \cg \nu ,\zeta \cd < 0 \p$$

Notons $\lambda \in \pic X$ le poids du faisceau $\call{L}$ (de sorte que $\call{L} \iso \call{L}_\lambda$).

\subsection{Cohomologie \`a support sur les vari\'et\'es stables}

Soit $I \sub\{1,...,r\}$ tel que $X_I= \adh{G.x}$. On note $X^+_I(x)$ la cellule $X^+(x) \cap X_I$ et $d_I$ la codimension de $X^+_I(x)$ dans $X_I$. 

On va d'abord montrer que le $\goth g$ module $H^{d_I}_{X_I^+(x)}(\call{L_\lambda}\res{X_I})$ n'a pas de sous-quotient simple de dimension finie en calculant son caractère, comme $\chap{T}-$module. 

\subsubsection{Calcul de caractères}

Avec les notations du lemme \ref{lem:tang}, on a :

\begin{pro}
Le $\chap{T}-$module $H^{d_I}_{X_I^+(x)}(\call{L_\lambda}\res{X_I})$ a pour caractère :
$$\frac{e^{\lambda'}}{\ma \prod_{i=1 \atop w(\delta_i) >0 }^u (1-e^{-w(\delta_i)}) \prod_{i=1 \atop w(\delta_i) <0 }^u (1-e^{w(\delta_i)}) \prod_{j=1}^v (1-e^{-w(\beta_j)})^2}$$
où $\ma \lambda':= \lambda_x - \sum_{i=1 \atop w(\delta_i) > 0}^u w (\delta_i) - \sum_{j=1}^v w(\beta_j)$ avec $\lambda_x$ le poids de la fibre $\call{L}_\lambda\res{x}$. 
\end{pro}

\begin{dem}
D'après \cite[prop. 11.9]{Kempf} (\cf aussi \cite[th. II.3.2]{these}), on a un isomorphisme de $\chap{T}-modules$ :
$$H^{d_I}_{X_I^+(x)}(\call{L}_\lambda\res{X_I}) \iso \Sym((T_x X_I)_+^*) \tenso{\kk} \Sym ((T_x X_I)_-) \tenso{\kk} \Vect^{d_I} (T_x X_I)_-$$
(où $\Sym$ est l'algèbre sym\'etrique et $\Vect^d$ la $d-$ième puissance ext\'erieure).

Or, d'après le lemme \ref{lem:tang}, les poids de $T_x X_I$ sont :$$ -w (\delta_1) ,..., -w(\delta_u), \pm w(\beta_1),..., \pm w(\beta_v) $$
et en particulier, ce sont des racines.
Parmi ces poids, ceux de $(T_x X_I)_+$ (resp. $(T_x X_I)_-$) sont les racines positives (resp. n\'egatives) car le sous-groupe \`a un paramètre $\zeta$ est dominant. De plus, comme $w \in W^{L_I}$ et comme les $\beta_j$ sont des racines positives de $L_I$, les poids $w(\beta_j)$ sont des racines positives.   
\end{dem}

{\bf Remarque :} Puisque le $\goth g-$module $H^{d_I}_{X_I^+(x)}(\call{L}_\lambda)$ a un caractère et puisque c'est aussi un $\chap{B}-$module (car la cellule $X_I^+(x)$ est $B-$invariante), c'est un $\goth g-$module de longueur finie. 

Nous allons voir que le caractère suffit pour d\'eterminer la multiplicit\'e des $\goth g-$modules simples de dimension finie dans le $\goth g-$module $H^{d_I}_{X_I^+(x)}(\call{L}_\lambda)$.

\subsubsection{Calcul de multiplicit\'es \`a partir de caractères}

On note $U(\cdot)$ l'algèbre enveloppante d'une algèbre de Lie.  On notera $Z(\goth g)$ le centre de l'algèbre $U(\goth g)$. Soit $\goth n^+$ l'algèbre de Lie de $B$.

On dit qu'un $\goth g-$module $M$ est $U(\goth n^+)-$fini si pour tout $m \in M$, le sous-espace $U(\goth n^+)m$ est de dimension finie. On dira que $M$ est $\chap{T}-$diagonalisable s'il est la somme directe de tous ses $\chap{T}-$espaces propres.

{\bf Exemple :} les $\goth g-\chap{B}-$modules sont $U(\goth n^+)-$finis et $\chap{T}-$diagonalisables.
\vskip .5cm

Le lemme qui suit montre en particulier comment v\'erifier la propri\'et\'e d'avoir un sous-quotient simple de dimension finie \`a partir du caractère. 

On notera, pour tout caractère $\lambda \in {\cal X}$, $L(\lambda)$ le $\goth g-$module simple de plus haut poids $\lambda$ et $\chi_\lambda : Z(\goth g) \to \kk$ son caractère central (on remarque que $L(\lambda)$ est de dimension finie si et seulement si $\lambda$ est dominant).

\begin{lem}\label{lem:car}
Soit $M$ un $\goth g-$module $U(\goth n^+)-$fini, $\chap{T}-$diagonalisable et avec un caractère $[M]$ dans l'anneau $\Z \cg {\cal X} \cd$. 

Alors pour tout caractère $\mu \in {\cal X}$ dominant, la multiplicit\'e du $\goth g-$module simple de plus haut poids $\mu$ est donn\'ee par :
$$[M : L(\mu)] = \prod_{\alpha \in \Phi^+} (1 - e^{-\alpha}) . [M] (\mu)$$
(la multiplicit\'e du poids $\mu$ dans le caractère virtuel $\ma \prod_{\alpha \in \Phi^+} (1 - e^{-\alpha}) . [M]$).
\end{lem}

\begin{dem}
Pour chaque caractère central $\chi : Z(\goth g) \to \K$, on note $M_\chi$ le sous-espace propre g\'en\'eralis\'e associ\'e \`a $\chi$ :

$$M_\chi:=\left\{ m \in M \tq \exist k >0 \tq  (\ker \chi)^k.m =0 \right\} \p$$

 Comme $M$ a un caractère, on a la d\'ecomposition :
$$M = \plus_{\chi} M_\chi$$
(\cf par exemple \cite[\S 8, pro. 8.6]{BGG}).

Il suffit donc de d\'emontrer le lemme lorsque $M = M_\chi$. Dans ce cas, $M = M_\chi$ admet une suite de Jordan-H\"{o}lder finie :
$$M=M_0 \con ... \con M_N \con M_{N+1} = 0$$
dont les sous-quotients simples sont tous de la forme :
$$M_i/M_{i+1} \iso L(\lambda^{(i)})$$
pour un certain $\lambda^{(i)} \in {\cal X}$ tel que $\chi_{\lambda^{(i)}} = \chi$ (\cf \cite[pro. 7.6.1]{Dix}).

Notons $\{\lambda_1,..., \lambda_t\}$ l'ensemble fini $\{\lambda \in {\cal X} \tq \chi_\lambda = \chi\}$.

  Pour chaque $i$, soit $n_i$ la multiplicit\'e du $\goth g-$module simple $L(\lambda_i)$ dans $M$. On a alors l'\'egalit\'e de caractères :\begin{equation}\label{eq:ML}
[M] = \sum_{1 \le i \le t} n_i [L(\lambda_i)] \p
\end{equation}

Or, les caractères des modules simples \`a plus haut poids s'expriment avec les caractères des modules de Verma : pour tout $i$, il existe des entiers $m_{i,j}$ (\'eventuellement n\'egatifs) tels que :
\begin{equation}\label{eq:tri}
[L(\lambda_i)] = [M(\lambda_i)] + \sum_{1 \le j \le t \atop \lambda_j < \lambda_i } m_{i,j} [M(\lambda_j)] \p 
\end{equation}

L'ensemble $\{\lambda_1,...,\lambda_t\}$ contient au plus un poids dominant ; supposons par exemple que $\lambda_1$ est l'unique poids dominant de cet ensemble. Dans ce cas, $\lambda_j < \lambda_1$ pour tout $2 \le j \le t$ et il r\'esulte de (\ref{eq:ML}) et (\ref{eq:tri}) que :
$$[M] = n_1 [M(\lambda_1)] + \sum_{2 \le j \le t} p_j [M(\lambda_j)]$$
pour certains entiers $p_j$. 

Connaissant les caractères des modules de Verma, on trouve donc :
$$[M] = \frac{n_1 e^{\lambda_1}}{\prod_{\alpha \in \Phi^+} (1-e^{-\alpha})} + \frac{p_2 e^{\lambda_2}}{\prod_{\alpha \in \Phi^+} (1-e^{-\alpha})} + ... +\frac{p_t e^{\lambda_t}}{\prod_{\alpha \in \Phi^+} (1-e^{-\alpha})} $$
\cad : \begin{equation}\label{eq:simpli}
\prod_{\alpha \in \Phi^+} (1-e^{-\alpha}) . [M] = n_1 e^{\lambda_1} + p_2 e^{\lambda_2} + ... + p_t e^{\lambda_t} \p
\end{equation}
En cons\'equence, $n_1$ est la multilicit\'e de $\lambda_1$ dans le caractère virtuel : $\prod_{\alpha \in \Phi^+} (1-e^{-\alpha}) [M]$.

Soit maintenant $\mu \in {\cal X}$ un caractère dominant. 

Si $\chi_\mu = \chi$, alors $\mu \in \{\lambda_1,...,\lambda_t\}$ et donc $\mu = \lambda_1$, d'où : $[M:L(\mu)] = n_1$.  

Si $\chi_\mu \not= \chi$, alors d'une part $[M:L(\mu)] = 0$ car $M = M_\chi$ et d'autre part : $\mu \not\in \{\lambda_1,...,\lambda_t\}$, ce qui entra\^ine :
$$\prod_{\alpha \in \Phi^+} (1-e^{-\alpha}) . [M] (\mu) = 0$$

d'après (\ref{eq:simpli}).
\end{dem}

\subsubsection{Absence de sous-quotients-simples de dimension finie}

Nous allons appliquer le lemme \ref{lem:car} au $\goth g-$module $ H^{d_I}_{X^+_I(x)}(\call{L_\lambda}\res{X_I})$ : 
\begin{pro}\label{pro:orbf}
Si $x$ n'est pas dans la $G-$orbite ferm\'ee de $X$, alors le $\goth g-$module $ H^{d_I}_{X^+_I(x)}(\call{L_\lambda})$ n'a pas de sous-quotient simples de dimension finie.
\end{pro}

\begin{dem}

Posons  $M :=H^{d_I}_{X^+_I(x)}(\call{L_\lambda}\res{X_I})$. Si $M$ a un sous-quotient simple de dimension finie $L(\mu)$ avec $\mu \in {\cal X}$, alors forc\'ement $\mu$ est dominant, et $\mu$ est un poids du caractère virtuel :
$$\prod_{\alpha \in \Phi^+} (1-e^{-\alpha}) [M] = \frac{(\prod_{\alpha \in \Phi^+} (1-e^{-\alpha})) e^{\lambda'}} {\ma \prod_{1\le i \le u \atop w\delta_i >0} (1-e^{-w\delta_i}) \prod_{1\le i \le u \atop w\delta_i <0} (1-e^{w\delta_i}) \prod_{1 \le j \le v} (1-e^{-w\beta_j})^2} \p$$

{ \bf Nous allons montrer que ce caractère virtuel n'a pas de poids dominant si $ v >0$}.

Notons $\beta_1,...,\beta_t$, avec $t \ge v$ les racines positives du sous-groupe de Levi $L_I$. On a la d\'ecomposition suivante :
$$\Phi^+= \{w\delta_i \tq 1 \le i \le u , w\delta_i >0\} \dij \{-w \delta_i \tq 1 \le i \le u, w\delta_i <0\} \dij \{w\beta_j \tq 1 \le j \le t\} \p$$

On a donc :
$$  \prod_{\alpha \in \Phi^+} (1-e^{-\alpha}) [M] = \frac{\ma (\prod_{j= v+1}^t (1-e^{-w\beta_j})) e^{\lambda'} }{\prod_{j=1}^v (1-e^{-w\beta_j})} \p$$

Par cons\'equent les poids du caractère $ \prod_{\alpha \in \Phi^+} (1-e^{-\alpha}) [M]$ sont de la forme $$\nu = \lambda'  - \sum_{j=1}^t n_j w\beta_j$$
$$
= \lambda_x - \sum_{i=1 \atop w \delta_i>0}^u w\delta_i - \sum_{j=1}^v w\beta_j - \sum_{j=1}^t n_j w \beta_j
$$
pour certains entiers $n_j \ge 0$ (rappelons que $\lambda_x$ d\'esigne le poids de la fibre $\call{L}_\lambda\res{x}$).

Posons $\gamma := \sum_{j=1}^v \beta_j + \sum_{j=1}^t n_j \beta_j$, de sorte que :
\begin{equation}\label{eq:formp}
\nu =\lambda_x - \sum_{i=1 \atop w \delta_i>0}^u w\delta_i - w\gamma
\end{equation}.

Comme $w \in W^{L_I}$, $w\gamma$ est une somme de racines positives. Donc, si on suppose que $\nu$ est un poids dominant, alors $(\nu,w\gamma) \ge 0$ ; cela signifie d'après (\ref{eq:formp}) que :

\begin{equation}\label{eq:nu}
 (\nu,w\gamma) = (\lambda_x, w\gamma) - (\sum_{i=1\atop w\delta_i >0}^u w\delta_i,w\gamma) -(w\gamma,w\gamma)  \ge 0 \p
\end{equation}

Mais, d'une part, on a :

\begin{equation}\label{eq:lwg}
(\lambda_x, w\gamma) =0 \p
\end{equation}
En effet, d'après le lemme \ref{lem:pfibre}, $\lambda_x$ est un caractère du groupe $wP_Iw\inv$ et $\gamma$ est une combinaison lin\'eaire de racines de $L_I$.

D'autre part, on a aussi :\begin{equation}\label{eq:srru}
(\sum_{i=1 \atop w\delta_i >0}^u w\delta_i, w \gamma ) \ge 0 \p
\end{equation}

En effet, les $\delta_i$ sont les racines de $R_u(P_I)$, \cad les racines positives qui ne sont pas racines de $(L_I,T)$, et donc si on note  $\rho_I$ la demi-somme des racines positives de $L_I$, on trouve  : 
$$
\sum_{i=1 \atop w\delta_i >0}^u w\delta_i = \sum_{\alpha \in \Phi^+ \moins \Phi^+_{L_I} \atop w \alpha >0} w\alpha $$
$$= \sum_{\beta > 0 \atop w\inv \beta \in \Phi^+ \moins \Phi^+_{L_I}} \beta$$
$$= \sum_{\beta > 0 \atop w\inv \beta \in \Phi^+} \beta - \sum_{\beta > 0 \atop w\inv \beta \in \Phi^+_{L_I}} \beta $$
$$=   \sum_{\beta > 0 \atop w\inv \beta \in \Phi^+} \beta - \sum_{\beta \in w \Phi^+_{L_I}} \beta$$
$$ = \sum_{\beta > 0 \atop w\inv \beta \in \Phi^+} \beta -2w\rho_I \p$$

Or, on a :
$$\rho + w\rho = \demi \sum_{\beta >0}\beta + \demi \sum_{\beta >0} w\beta$$
$$= \demi\sum_{\beta >0 \atop w\inv \beta >0} \beta + \demi\sum_{\beta >0 \atop w\inv \beta <0} \beta + \demi\sum_{\beta  \atop w\inv \beta >0} \beta$$
$$= \sum_{\beta >0 \atop w\inv \beta >0} \beta + \demi\sum_{\beta >0 \atop w\inv \beta <0} \beta + \demi\sum_{\beta < 0  \atop w\inv \beta >0} \beta $$
$$= \sum_{\beta >0 \atop w\inv \beta >0} \beta \p$$
Il en r\'esulte que :
$$\sum_{i=1 \atop w\delta_i >0}^u w\delta_i = \rho + w\rho-2w\rho_I \p$$

Mais alors, on obtient :
$$(\sum_{i=1 \atop w\delta_i >0}^u w\delta_i , w\gamma ) = (\rho + w\rho-2w\rho_I , w \gamma) \p$$

Si $\beta_j$ est une racine simple de $L_I$, alors on a :
$$ \cg \rho + w\rho-2w\rho_I , w \beta_j^\ch \cd =  \cg \rho , w \beta_j^\ch  \cd + \cg \rho ,\beta_j^\ch  \cd  -2  \cg \rho_I ,\beta_j^\ch  \cd $$
$$ = \cg \rho , w \beta_j^\ch  \cd -1 \ge 0$$
car $w\beta_j$ est une racine positive.

Donc : 
$$(\rho + w\rho-2w\rho_I , w \beta_j )\ge 0$$
pour toute racine simple $\beta_j$ de $L_I$ et, par cons\'equent, aussi pour toute racine positive $\beta_j$ de $L_I$.

Et finalement :
$$(\sum_{i=1 \atop w\delta_i >0}^u w\delta_i , w\gamma ) = (\rho + w\rho-2w\rho_I , w \gamma)$$
$$= \sum_{j=1}^v (\rho + w\rho-2w\rho_I , w\beta_j) + \sum_{j=1}^t n_j (\rho + w\rho-2w\rho_L , w\beta_j) \ge 0$$
d'où (\ref{eq:srru}).

On d\'eduit de (\ref{eq:nu}), (\ref{eq:lwg}) et (\ref{eq:srru}) que si $\nu$ est un poids dominant, alors :

$$(w\gamma , w\gamma) \le 0$$
\ie : $\gamma = 0$. Comme $\gamma =\sum_{j=1}^v \beta_j + \sum_{j=1}^t n_j \beta_j$, on a n\'ecessairement $v =0$.

\begin{center}
***
\end{center}

En conclusion, puisque $2v$ est la codimension du groupe d'isotropie $G_x$ dans le sous-groupe parabolique $wP_I w\inv$ (\cf le lemme \ref{lem:tang}), on a : $G_x = wP_Iw\inv$ et finalement $G.x$ est projective donc ferm\'ee ; d'où la proposition \ref{pro:orbf}. 
\end{dem}

\subsection{Filtrations}
   
Pour terminer la d\'emonstration du th\'eorème principal \ref{th:orbf}, on utilise le r\'esultat suivant qui nous ramène \`a l'\'etude, faite ci-dessus, du $\goth g-$module $H^{d_I}_{X_I^+ (x)}(\call{L}_\lambda\res{X_I})$ :

\begin{thm}[{\cite[th. 4.1]{toi}}]
Le $\goth g-$module $H^d_{X^+(x)}(\call{L}_\lambda)$ est de longueur finie et possède une filtration de sous-$\goth g-$modules :
$$H^d_{X^+(x)}(\call{L}_\lambda) = M_0 \con M_1 \con ... \con M_n \con ...$$
telle que $\bigcap_{k \ge 0} M_k = 0$ et pour tout $k \ge 0$ :
$$ M_k/M_{k+1} = H^{d_I}_{X^+_I(x)} (\call{L}_\lambda \tens \call{O}_X(D_k)\res{X_I})$$
 où $D_k \in \sum_{i=1}^r \Z D_\alpha$ est un diviseur de $X$.
\end{thm}

Puisque les faisceaux inversibles $\call{L}_\lambda \tens\call{O}_X(D_k)$ sont encore sp\'eciaux, la proposition \ref{pro:orbf} s'applique : on en d\'eduit que les quotients successifs de la filtration du th\'eorème ci-dessus n'ont aucun  sous-quotient simple de dimension finie lorsque le point $x$ est hors de l'orbite ferm\'ee. Il en est donc de m\^eme pour le $  \goth g-$module $H^d_{X^+(x)}(\call{L}_\lambda)$.

Cela achève la d\'emonstration du th\'eorème \ref{th:orbf}.

\section{Applications}

Nous allons appliquer le th\'eorème \ref{th:orbf} au problème du calcul des groupes de cohomologie $H^d(X,\call{L}_\lambda)$, $d \ge 0, \lambda \in \pic(X)$.

Comme les faisceaux $\call{L}_\lambda$ sont $\chap{G}-$lin\'earis\'es sur $X$, tous ces groupes de cohomologie sont des $\chap{G}-$modules. Donc pour chaque entier $d$ et chaque poids $\lambda \in \pic X$, on a une d\'ecomposition en somme directe de $\chap{G}-$modules simples :
$$H^d(X,\call{L}_\lambda) = \Plus_{\mu \in {\cal X} \atop \mu \mathrm{\;dominant}} L(\mu)^{\plus m^d_\lambda(\mu)}$$
où $m^d_\lambda(\mu) \in \Z_{\ge 0}$ est la multiplicit\'e du module simple $L(\mu)$ dans le $\chap{G}-$module  $H^d (X,\call{L}_\lambda)$.

\subsection{Suite spectrale de Grothendieck-Cousin}\label{sec:ssgc}

Pour estimer les entiers $m^d_\lambda(\mu)$, on dispose d'une suite spectrale qui fait intervenir des groupes de cohomologie \`a support dans des cellules de Bialynicki-Birula de $X$.

Fixons une d\'ecomposition cellulaire de $X$ :
$$X = \Dij_{i=0}^N X^+_i$$
\cad que l'on note $x_0,...,x_N$ les points fixes de $T$ dans $X$, que l'on choisit un sous-groupe \`a un paramètre $\zeta$ de $T$ tel que $X^{\zeta(\kk^*)} = X^T$ et que l'on pose $X^+_i := X^+(x_i)$ la cellule de Bialynicki-Birula centr\'ee en $x_i$. On supposera de plus que $\zeta$ est dominant et r\'egulier \ie :

$$\qq \alpha \in \Phi^+, \cg \alpha, \zeta \cd >0 \p$$

D'après \cite[th. 3]{BB76}, il existe une suite d\'ecroissante de sous-espaces ferm\'es de $X$ :
$$X=Z_0 \con ... \con Z_N \con Z_{N+1}= \vide$$
telle que pour chaque $i$, $Z_i \moins Z_{i+1} = X^+_i$ (quitte \`a renum\'eroter les points fixes $x_i$).

D'après \cite{Kempf}, on a alors une suite spectrale de $\goth g-$modules qui converge :
\begin{equation}\label{eq:suispe}
E_1^{p,q} = H^{p+q}_{X^+_p}(\call{L}_\lambda) \impliq H^{p+q}( X, \call{L}_\lambda) \p \end{equation}

Les termes initiaux sont seulement des $\goth g-$modules mais les termes finaux sont des $\chap{G}-$modules, \cad des sommes directes de $\goth g-$modules de dimension finie.

A priori cette suite spectrale ne d\'eg\'enère pas mais on peut la simplifier. en effet, tous les termes $E_1^{p,q}$ tels que $p+q \not=\codim_X (X^+_p)$ sont nuls. D'un autre c\^ot\'e, gr\^ace au th\'eorème \ref{th:orbf}, une autre simplification est possible : on va pouvoir ne tenir compte que des $E_1^{p,q}$ tels que le point fixe $x_p$ est dans l'orbite ferm\'ee. En effet, les autres termes $E_1^{p,q}$ n'ont pas de $\chap{G}-$modules parmi leurs sous-quotients simples.

En ce qui concerne les caractères de $\chap{T}$, on utilisera les notations suivantes :
pour tout $\mu \in {\cal X}$, il existe au plus un $w \in W$ tel que $w(\mu + \rho) - \rho$ soit dominant ; dans ce cas, on note :
$$\mu^+:= w (\mu + \rho) - \rho \et l(\mu) := l(w)=|\left\{\alpha \in \Phi^+ \tq \cg \mu+\rho , \alpha^\ch \cd < 0\right\}|\p $$ 

\subsection{Caract\'eristique d'Euler-Poincar\'e}

Si $V$ est un $\chap{G}-$module rationnel, on note $\cl(V)$ la classe des $\chap{G}-$modules isomorphes \`a $V$. Le {\it groupe de Grothendieck} des $\chap{G}-$modules rationnels, not\'e $K(\chap{G})$, est le groupe commutatif d\'efini par les g\'en\'erateurs : $\cl(V)$ et par les relations  : $\cl(V) = \cl (V') + \cl (V'')$ pour chaque suite exacte courte $0 \to V' \to V \to V''$ de $\chap{G}-$modules. Pour tout $\chap{G}-$module $V$, on note $[V]$ l'image de $V$ dans $K(\chap{G})$.

Pour tout $\lambda\in \pic X$, notons $\chi(X,\call{L}_\lambda$) l'\'el\'ement
$$\sum_{d \ge 0} (-1)^d [H^d(X, \call{L}_\lambda)]$$ de $K(\chap{G})$ ; c'est la caract\'eristique d'Euler-Poincar\'e globale du faisceau inversible $\call{L}_\lambda$. 

Pour tout $\lambda\in \pic X$ et tout $\mu \in {\cal X}$ dominant, notons :
$$\chi_\lambda(\mu):= \sum_{d \ge 0} (-1)^d m^d_\lambda(\mu) \;\;,$$
c'est la multiplicit\'e selon $[L(\mu)]$ du $\chap{G}-$module virtuel $\chi(X,\call{L}_\lambda)$.

\begin{thm}\label{thm:euler}
Pour tout $1 \le i \le l$, soit $\rho_i := \frac{(\rho , \tilda{\alpha_i})}{(\tilda{\omega_i},\tilda{\alpha_i})}$. Pour toute partie $J$ de $\{1,...,r\}$, soient :
$$R_J:= \sum_{i=1 \atop i \in J}^r \Z_{>0} \tilda{\alpha_i} + \sum_{i=1 \atop i \not\in J}^r \Z_{\le 0} \tilda{\alpha_i} $$
$$\Omega_J := \sum_{i=1 \atop i \in J}^r \Z_{<-\rho_i} \tilda{\omega_i} + \sum_{i=1 \atop i \not\in J}^r \Z_{> - \rho_i} \tilda{\omega_i} \p$$

Si $\lambda$ est un poids sp\'ecial de $\pic(X)$, alors on a :
$$\chi(X,\call{L}_\lambda) = \sum_{J \sub \{1,...,r\}} \sum_{\nu \atop \nu + \rho \mbox{ \scriptsize r\'egulier}}(-1)^{l(\nu) +|J|} [L(\nu^+)]$$
où dans la deuxième somme, $\nu$ d\'ecrit l'ensemble de caractères $\lambda + R_J \cap \Omega_J$.

Autrement dit, pour tout poids dominant $\mu \in {\cal X}$ :
$$\chi_\lambda(\mu) = \sum_{J \sub \{1,...,r\}} \sum_{\nu} (-1)^{l(\nu) +|J|} $$
où cette fois, $\nu$ d\'ecrit l'ensemble fini de caractères $\lambda + R_J \cap \Omega_J \cap W *\mu$.
\end{thm}

\begin{figure}[h]
\begin{center}
\input{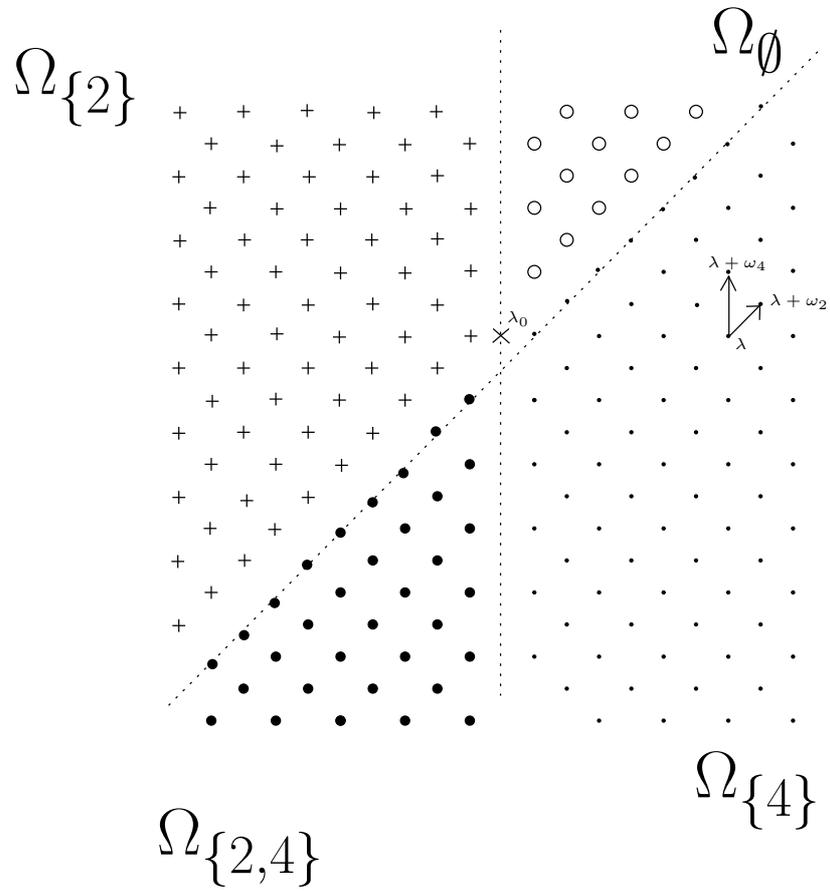}
    \caption{les ensembles $\Omega_J$ pour la compactification magnifique de $\Sp_{2n}/\Sp_{4} \croi \Sp_{2n-4}$, $n \ge 4$ où : $\lambda_0 := \frac{(\rho,\tilda{\alpha_2})}{(\omega_2,\tilda{\alpha_2})} \omega_2 +\frac{(\rho,\tilda{\alpha_4})}{(\omega_4,\tilda{\alpha_4})}\omega_4 = 2 \omega_2 + \frac{2n-5}{2}\omega_4$}
    \label{fig:c2P}
\end{center}
  \end{figure}

\begin{figure}[h]
\begin{center}
\input{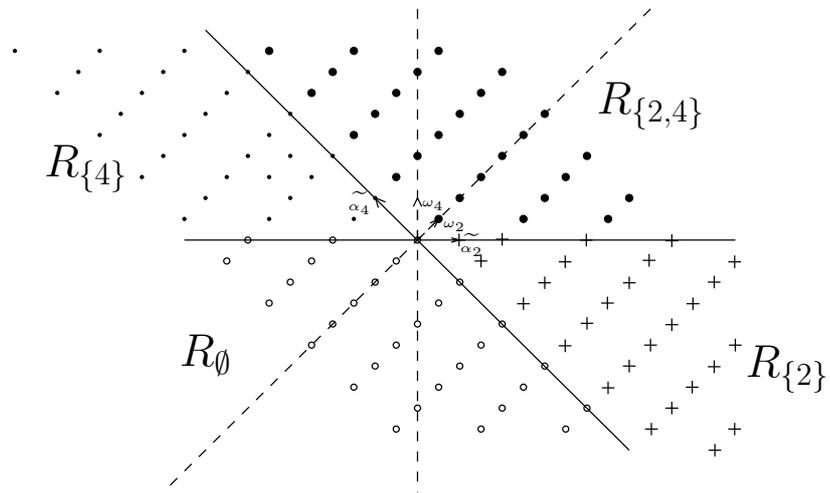}
    \caption{ les ensembles $R_J$ pour la compactification magnifique de $\Sp_{8}/\Sp_{4} \croi \Sp_{4}$}
    \label{fig:c2R4}
\end{center}
  \end{figure}

{\bf Remarque :} Sur les figures \ref{fig:c2P},\ref{fig:c2R} et \ref{fig:c2R4}, on a repr\'esent\'e les ensembles $\Omega_J$ et $R_J$ dans le cas de la compactification magnifique de l'espace homogène $\Sp_{2n}/\Sp_4 \croi \Sp_{2n-4}$, $n \ge 4$ (\cf \S \ref{sec:ega} pour les notations). 

\begin{dem}
Pour calculer la caract\'eristique d'Euler-Poincar\'e globale des faisceaux $\call{L}_\lambda$, on a seulement besoin des premiers termes de la suite spectrale (\ref{eq:suispe}).

On veut calculer :
$$\chi_\lambda(\mu) = \sum_{d \ge 0} (-1)^d[H^d(X,\call{L}_\lambda) : L(\mu)]\p$$

D'après (\ref{eq:suispe}), on a :
$$ \chi_\lambda(\mu) = \sum_{p,q} (-1)^{p+q} [E_1^{p,q} : L(\mu)]$$
$$= \sum_{d \ge 0} (-1)^d \sum_{p \ge 0} [H^d_{X^+_p}(\call{L}_\lambda):L(\mu)] $$

(rappelons que les cellules $X^+_p$ sont d\'efinies par un sous-groupe \`a un paramètre $\zeta$ dominant et r\'egulier, \cf \S \ref{sec:ssgc}).

Or, d'après le th\'eorème principal \ref{th:orbf}, si la cellule $X^+_p$ n'est pas centr\'ee en un point de l'orbite ferm\'ee, alors la multiplicit\'e $\ma [H^d_{X^+_p}(\call{L})_\lambda):L(\mu)]$ est nulle. Donc, si l'on note $X^+_w$ la cellule de Bialynicki-Birula centr\'ee en le point $wQ/Q$, on trouve :
$$\chi_\lambda(\mu) = \sum_{d \ge 0} (-1)^d \sum_{w \in W^Q} [H^d_{X^+_w}(\call{L}_\lambda):L(\mu)] $$

(Rappelons que les points fixes de l'orbite ferm\'ee $F = G/Q$ sont param\'etr\'es par $W^Q$, l'ensemble des repr\'esentants de longueur minimale de $W/W_Q$).

D'un autre c\^ot\'e, le groupe de cohomologie \`a support $H^d_{X^+_w}(\call{L}_\lambda)$ est nul si $d \not= \codim_X X^+_w$ et si $d = \codim_X X^+_w$, on a d'après \cite[pro. 4.6 et th. 4.4]{toi} :\begin{equation}\label{eq:multcs}
[H^d_{X^+_w}(\call{L})_\lambda):L(\mu)]
\end{equation}
$$= \left\{ \begin{array}{cl}
1 & \ma \si w\inv(\mu + \rho) \in \lambda + \rho + \sum_{i=1\atop X^+_w \sub D_i }^r \Z_{> 0} \tilda{\alpha}_i + \sum_{i=1\atop X^+_w \not\sub D_i }^r \Z_{\le  0} \tilda{\alpha}_i \\
0 & \mbox{ sinon.}
\end{array}\right.$$

Or, pour chaque cellule $X^+_w$ et chaque diviseur $D_i$, le tore $\chap{T}$ agit sur la fibre ${\cal O}_X(D_i)\res{wQ/Q}$ via le caractère $w(\tilda{\alpha_i})$. On obtient donc les \'equivalences suivantes :
$$X^+_w \sub D_i \equi \cg w(\tilda{\alpha_i}) , \zeta \cd < 0$$
$$X^+_w \not\sub D_i \equi \cg w(\tilda{\alpha_i}) , \zeta \cd > 0 \p$$
De plus, comme $X^+_w \cap \bigcap_{i=1}^r D_i = BwQ/Q$, on a :
\begin{equation}\label{eq:codimcell}
\codim_X X^+_w = \left | \{ 1 \le i \le r \tq \cg w(\tilda{\alpha_i}) , \zeta \cd < 0 \} \right | + l(w)   \p
\end{equation}
Posons $J_w := \{ 1 \le i \le r \tq \cg w(\tilda{\alpha_i}) , \zeta \cd < 0 \} + l(w) \}$.
On d\'eduit de ce qui pr\'ecède que :
$$\chi_\lambda(\mu) = \sum_w (-1)^{l(w) + |J_w|} $$
où $w$ d\'ecrit l'ensemble, que nous noterons $W_{\lambda,\mu}$, des \'el\'ements $w$ de $W$ tels que :
$$w\inv(\mu + \rho) \in \lambda + \rho + \sum_{i=1\atop i \in J_w }^r \Z_{> 0} \tilda{\alpha}_i + \sum_{i=1 \atop i \not\in J_w }^r \Z_{\le  0} \tilda{\alpha}_i \p$$

Si on admet le lemme suivant :
\begin{lem}\label{lem:wspecial}
Pour tout poids sp\'ecial $\lambda$ tel que $\lambda + \rho$ est r\'egulier et pour tout $\alpha \in \Delta_1$ :
$$\cg w_\lambda (\tilda{\alpha}) , \zeta \cd < 0 \equi (\lambda + \rho , \tilda{\alpha}) < 0 \p$$ 
\end{lem}

alors on remarque que l'application
$$
\appli{W_{\lambda,\mu}}{\{(J,\nu) \tq J \sub \{1,...,r\} ,\;\nu \in \lambda + R_J \cap \Omega_J \cap W*\mu \}}{w}{(J,w\inv (\mu+\rho) - \rho)}
$$
est une bijection. D'où les formules de l'\'enonc\'e.

D\'emontrons maintenant le lemme \ref{lem:wspecial} :

L'\'equivalence \`a d\'emontrer n'est pas imm\'ediate car a priori, $w_\lambda(\tilda{\alpha})$ n'est pas un \'el\'ement de $\tilda{\Phi_1}$. N\'eanmoins, on va montrer que $w_\lambda(\tilda{\alpha})$ est une somme de $2$ racines positives ou de $2$ racines n\'egatives. 

Pour cela, on v\'erifie qu'il existe des entiers $n_\delta \ge 0, \delta \in \Delta_0$ tels que :
\begin{equation}\label{eq:mqu}
\alpha':=\alpha +\sum_{\delta \in \Delta_0} n_\delta \delta \in \Phi_1^+ \et
|\cg \rho , \alpha'^\ch \cd + \cg \rho , \theta(\alpha')^\ch \cd | \le 1 \p
\end{equation} 

Il suffit, en effet, de traiter le cas où $X$ est de rang $1$ avec $\tilda{\Phi_1} = \{\pm \tilda{\alpha}\}$. Cela fait seulement huit possibilit\'es ; dans le tableau suivant, chaque vari\'et\'e sym\'etrique complète est repr\'esent\'ee par le diagramme de Satake de l'espace homogène ${G}/{H}$ (c'est le diagramme de Dynkyn de $G$ où les sommets correspondant aux \'el\'ements de $\Delta_0$ sont noircis et où les racines simples $\alpha_i$ et $-\adh{\theta}(\alpha_i)$ distinctes sont reli\'ees par \flfl). Pour chaque exemple, on donne $\theta(\alpha)$ et une racine $\alpha'$ qui v\'erifie (\ref{eq:mqu}) :
 $$
\begin{array}{|c|c|}
\hline
X & \theta \atop \alpha'\\
\hline
\stackrel{\alpha_1}{\circ} & -\theta(\alpha_1) = \alpha_1 \atop \alpha' =  \alpha_1 \\
\hline
\xymatrix{\stackrel{\alpha_1}{\circ}\ar@/^/ @{<->}[r] &\stackrel{\alpha_2}{\circ}} & \theta(\alpha_1) = \alpha_2 \atop \alpha' =  \alpha_1\\
\hline
\xymatrix{\stackrel{\alpha_1}{\bullet} \ar@{-}[r]<-.7ex> & \stackrel{\alpha_2}{\circ} \ar@{-}[r] <-.7ex> & \stackrel{\alpha_3}{\bullet}} & -\theta(\alpha_2)= \alpha_1+\alpha_2+\alpha_3 \atop \alpha' =  \alpha_1 + \alpha_2\\
\hline
\xymatrix{\stackrel{\alpha_1}{\circ} \ar@{-}[r] <-.7ex> & \stackrel{\alpha_2}{\bullet} \ar@{.}[r] <-.7ex> &\stackrel{\alpha_{n-1}}{\bullet} \ar@{=>}[r] <-.7ex> & \stackrel{\alpha_n}{\bullet}} \; n \ge 2& -\theta(\alpha_1) = \alpha_1 + 2\sum_{i=2}^n \alpha_i \atop \alpha' =  \alpha_1 +  \sum_{i=2}^n \alpha_i\\
\hline
\xymatrix{& & & \stackrel{\alpha_{n-1}}{\bullet} \\
\stackrel{\alpha_1}{\circ} \ar@{-}[r] <-.7ex> & \stackrel{\alpha_2}{\bullet} \ar@{.}[r] <-.7ex> & \stackrel{\alpha_{n-2}}{\bullet} \ar@{-}[ur] \ar@{-}[dr] & \;\; n \ge 4\\
& & & \stackrel{\alpha_{n}}{\bullet} } & -\theta(\alpha_1) = \alpha_1 + 2\sum_{i=2}^{n-2} \alpha_i + \alpha_{n-1} + \alpha_n \atop \alpha' =  \sum_{i=1}^n \alpha_i \\
\hline
\xymatrix{\stackrel{\alpha_1}{\circ}\ar@{-}[r] <-.7ex> \ar@/^2ex/ @{<->}[rrr]&  \stackrel{\alpha_2}{\bullet} \ar@{.}[r] <-.7ex> & \stackrel{\alpha_{n-2}}{\bullet}  \ar@{-}[r] <-.7ex> & \stackrel{\alpha_n}{\circ} } \; n \ge 3& -\theta(\alpha_1) = \sum_{i=2}^n \alpha_i \atop \alpha' =  \sum_{i=1}^{[\frac{n}{2}]} \alpha_i\\
\hline
\xymatrix{\stackrel{\alpha_1}{\bullet} \ar@{-}[r] <-.7ex> &  \stackrel{\alpha_2}{\circ} \ar@{-}[r] <-.7ex> & \stackrel{\alpha_{3}}{\bullet}  \ar@{.}[r] <-.7ex> & \stackrel{\alpha_{n-1}}{\bullet} \ar@{<=}[r] <-.7ex> & \stackrel{\alpha_{n}}{\bullet} } \; n \ge 2& -\theta(\alpha_2) = \alpha_1 + \alpha_2 + 2 \sum_{i=3}^{n-1} \alpha_i + \alpha_n \atop \alpha' =  \sum_{i=1}^{n-1} \alpha_i\\
\hline
\xymatrix{\stackrel{\alpha_1}{\bullet} \ar@{-}[r] <-.7ex> & \stackrel{\alpha_2}{\bullet} \ar@{=>}[r] <-.7ex> &\stackrel{\alpha_3}{\bullet} \ar@{-}[r] <-.7ex> & \stackrel{\alpha_4}{\circ}  } & -\theta(\alpha_4) = \alpha_1+2\alpha_2+3\alpha_3 + \alpha_4 \atop \alpha' =  \alpha_1+ \alpha_2 + \alpha_3 + \alpha_4\\
\hline
\end{array}$$

On choisit donc $\alpha'$ qui v\'erifie (\ref{eq:mqu}) : on a alors $\tilda{\alpha} = \tilda{\alpha'}$ et $w_\lambda(\tilda{\alpha}) = w_\lambda(\alpha') - w_\lambda(\theta(\alpha'))$. 

Or,  comme le caractère $w_\lambda(\lambda+\rho)$ est dominant r\'egulier et comme $\theta(\lambda) = - \lambda$, on a les \'equivalences suivantes :

$$w_\lambda(\alpha') \in \Phi^+ \equi \cg w_\lambda(\lambda + \rho), w_\lambda(\alpha')^\ch \cd > 0 $$
$$ \equi \cg \lambda , \alpha'^\ch \cd + \cg \rho , \alpha'^\ch \cd >0 $$
et :
$$- w_\lambda(\theta(\alpha')) \in \Phi^+ \equi \cg \lambda + \rho , -\theta(\alpha')^\ch \cd >0$$
$$ \equi \cg \lambda , -\theta(\alpha')^\ch \cd + \cg \rho ,  -\theta(\alpha')^\ch \cd > 0$$
$$ \equi \cg \lambda , \alpha'^\ch \cd - \cg \rho ,  \theta(\alpha')^\ch \cd > 0 \p$$

Mais puisque, le poids $\lambda + \rho$ est entier et r\'egulier, les nombres 
$$\cg \lambda , \alpha'^\ch \cd + \cg \rho , \alpha'^\ch \cd \et  \cg \lambda , \alpha'^\ch \cd - \cg \rho , \theta(\alpha')^\ch \cd$$ 
sont des entiers non nuls, dont la diff\'erence est en valeur absolue inf\'erieure ou \'egale \`a $1$, d'après (\ref{eq:mqu}). Ce sont donc des entiers de m\^emes signes.

Comme le sous-groupe \`a un paramètre $\zeta$ est dominant r\'egulier, on peut conclure :
$$\cg w_\lambda(\tilda{\alpha}) , \zeta \cd >0 \equi w_\lambda(\alpha') \et w_\lambda(-\theta(\alpha')) \in \Phi^+$$
$$\equi \cg \lambda + \rho , \alpha'^\ch \cd \et \cg \lambda + \rho -\theta(\alpha')^\ch \cd >0 $$
$$ \equi (\lambda + \rho , \alpha') \et (\lambda+ \rho , -\theta(\alpha')) >0 $$
$$ \equi (\lambda + \rho, \tilda{\alpha'}) >0$$
$$ \equi (\lambda + \rho, \tilda{\alpha}) >0 \p$$
\end{dem}

{\bf Remarque :} En particulier, les codimensions des cellules de Bialynicki-Birula qui interviennent dans le calcul de la caract\'eristique d'Euler-Poincar\'e sont ind\'ependantes du sous-groupe \`a un paramètre $\zeta$ choisi.

\subsection{Estim\'ee des multiplicit\'es}
Pour deux $\chap{G}-$modules $V_1,V_2$ de dimension finie, la notation $V_1 \le V_2$ signifiera que pour tout caractère dominant $\mu$ :
$$[V_1:L(\mu)] \le [V_2:L(\mu)] \p$$

On obtient gr\^ace \`a la suite spectrale (\ref{eq:suispe}) et avec les notations du th\'eorème \ref{thm:euler} la majoration suivante des multiplicit\'es $m^d_\lambda(\mu)$ :
\begin{thm}\label{thm:ineg}
Pour tout poids sp\'ecial $\lambda \in \pic(X)$ et pour tout entier $d$ :

\begin{equation}\label{eq:ineg}
H^d(X,\call{L}_\lambda) \le \Plus_{J \sub \{1,...,r\}} \Plus_{{\nu \in \lambda + R_J \cap \Omega_J \atop \nu + \rho \mbox{ \scriptsize r\'egulier }} \atop l(\nu) +|J| = d} L(\nu^+) \p
\end{equation}

\end{thm}

{\bf Remarque :} Lorsque $X$ est de rang minimal, \ie lorsque $r = \rg(G) - \rg (H)$, cette in\'egalit\'e est une \'egalit\'e (\cf \cite[th. 3.1]{toi2}). Nous donnerons dans la section suivante un exemple, qui n'est pas de rang minimal, pour lequel on peut aussi d\'emontrer l'\'egalit\'e.  En revanche, nous verrons que cette in\'egalit\'e peut devenir stricte dans le cas de la vari\'et\'e des coniques complètes (\cf le th\'eorème \ref{thm:ccc}).

\begin{dem}
On d\'eduit de (\ref{eq:suispe}) que pour tout poids dominant $\mu$ :
$$[H^d(X,\call{L}_\lambda) : L(\mu)] \le \sum_{p,q \atop p+q = d}[E_1^{p,q} : L(\mu)]$$
$$ \le \sum_{p,q \atop p+q = d} [H^d_{X^+_p}(\call{L}_\lambda) : L(\mu)] \p$$

Or, d'après le th\'eorème \ref{th:orbf}, la multiplicit\'e  $[H^d_{X^+_p}(\call{L}_\lambda) : L(\mu)]$ est nulle si la cellule $X^+_p$ n'est pas centr\'ee en un point de l'orbite ferm\'ee $F \iso G/Q$. Donc :
$$[H^d(X,\call{L}_\lambda) : L(\mu)] \le \sum_{w \in W^Q} [H^d_{X^+_w}(\call{L}_\lambda) : L(\mu)] \p$$

L'in\'egalit\'e de l'\'enonc\'e r\'esulte alors de (\ref{eq:multcs}) et du lemme \ref{lem:wspecial}.
\end{dem}

On d\'eduit de ce th\'eorème le r\'esultat d'annulation suivant :
\begin{cor}\label{cor:annul1}
Pour tout poids sp\'ecial $\lambda \in \pic(X)$, $H^1(X,\call{L}_\lambda) =0$.
\end{cor}

{\bf Remarque : }En particulier, si $X$ est une vari\'et\'e sym\'etrique complète non exceptionnelle (\ie si tous les poids de $\pic X$ sont  sp\'eciaux ) alors, pour tout faisceau invercible $\call{L}$ sur $X$, on a : $H^1(X,\call{L}) = 0$. C'est le cas par exemple pour la {\it vari\'et\'e des quadriques complètes}, la compactification magnifique de $PGL_n/PSO_n$.

\begin{dem}[du corollaire]
Soient $J \sub \{1,...,r\}$ et $\nu \in \lambda + R_J \cap \Omega_J$ tel que $\nu + \rho$ est r\'egulier. 

Nous allons montrer que $l(\nu) + |J|$ ne vaut jamais $1$ :

Remarquons que :
$$\nu \in \lambda + R_J \impliq \theta(\nu) = -\nu \p$$
On en d\'eduit :
$$\qq \delta \in \Delta_0,\; \cg \nu , \delta^\ch \cd = \cg \theta(\nu), \theta(\delta)^\ch \cd = - \cg \nu, \delta^\ch \cd$$
\begin{equation}\label{eq:nuort}
 \impliq  \qq \delta \in \Delta_0,\; \cg \nu , \delta^\ch \cd = 0 \p
\end{equation}

On a :
$$J = \{1 \le i \le r \tq (\nu + \rho , \tilda{\alpha_i}) < 0 \} \et  l(\nu) = \left | \{\alpha \in \Phi^+ \tq (\nu + \rho , \alpha) < 0\} \right |\p$$

Donc : $|J| \le l(\nu)$. Mais alors :
$$l(\nu) + |J| = 1 \impliq |J| = 0 \et l(\nu) =1 \p$$

Nous allons voir que dans ce cas, on a en fait $l(\nu) \ge 2$ : absurde.

D'une part, il existe un unique $\alpha \in \Delta$ tel que $w_\nu(\alpha) \in \Phi^-$ \ie : $\cg \nu + \rho , \alpha^\ch \cd < 0$ ou encore :
\begin{equation}\label{eq:nupe}
\cg \nu , \alpha^\ch \cd + 1 <0 
\end{equation}
(en particulier, $\alpha \in \Delta_1$).

D'autre part, comme $|J| = 0$, on a : $(\nu + \rho , \tilda{\alpha} ) > 0$ et donc $ \cg \nu + \rho , -\theta(\alpha)^\ch \cd > 0$. Par cons\'equent :
$$\cg \rho , -\theta(\alpha)^\ch \cd > - \cg \nu , -\theta(\alpha)^\ch \cd$$
$$ \impliq \cg \rho , -\theta(\alpha)^\ch \cd > - \cg \nu , \alpha^\ch \cd > 1 \p$$

On est donc sous l'hypothèse du lemme suivant :

\begin{lem}\label{lem:rhoth}
Soit $\alpha \in \Delta_1$ tel que $\cg \rho , - \theta( \alpha)^\ch \cd > 2$. Alors il existe $\delta \in \Delta_0$ tel que :
$$\alpha + \delta \in \Phi^+, \; \cg \rho , (\alpha + \delta )^\ch \cd = 1 \mbox{ ou }2 , \; (\alpha+ \delta, \alpha + \delta) = (\alpha, \alpha) \p$$
\end{lem}

Avec (\ref{eq:nuort}) et (\ref{eq:nupe}), ce lemme permet de conclure :

$$\cg \nu + \rho , (\alpha + \delta)^\ch \cd \le \cg \nu , (\alpha + \delta)^\ch \cd + 2$$
$$\le  2 \frac{(\nu , \alpha + \delta)}{(\alpha + \delta , \alpha + \delta)} + 2$$
$$\le 2\frac{(\nu , \alpha )}{(\alpha + \delta , \alpha + \delta)} + 2$$
$$ \le \cg \nu, \alpha^\ch \cd + 2 \le 0 $$
$$\impliq \cg \nu + \rho , (\alpha + \delta)^\ch \cd < 0$$
car $\nu + \rho$ est r\'egulier ; et donc, $l(\nu) \ge 2$.

Il reste \`a d\'emontrer le lemme \ref{lem:rhoth} :

il suffit de traiter les $6$ cas où $X$ est de rang $1$ et où $\Delta_0 \not = \vide$. On le fait dans le tableau suivant en donnant pour chaque exemple, une racine $\alpha + \delta$ qui v\'erifie l'\'enonc\'e du lemme :
\begin{center}$$
\begin{array}{|c|c|c|c|}
\hline
X & \alpha & \alpha + \delta & \cg \rho , (\alpha + \delta)^\ch \cd \\
\hline
\xymatrix{\stackrel{\alpha_1}{\bullet} \ar@{-}[r]<-.7ex> & \stackrel{\alpha_2}{\circ} \ar@{-}[r] <-.7ex> & \stackrel{\alpha_3}{\bullet}} & \alpha_2 & \alpha_1 + \alpha_2 & 2\\
\hline
\xymatrix{\stackrel{\alpha_1}{\circ} \ar@{-}[r] <-.7ex> & \stackrel{\alpha_2}{\bullet} \ar@{.}[r] <-.7ex> &\stackrel{\alpha_{n-1}}{\bullet} \ar@{=>}[r] <-.7ex> & \stackrel{\alpha_n}{\bullet}} \; n \ge 3 & \alpha_1 & \alpha_1 + \alpha_2 & 2\\
\hline
\xymatrix{& & & \stackrel{\alpha_{n-1}}{\bullet} \\
\stackrel{\alpha_1}{\circ} \ar@{-}[r] <-.7ex> & \stackrel{\alpha_2}{\bullet} \ar@{.}[r] <-.7ex> & \stackrel{\alpha_{n-2}}{\bullet} \ar@{-}[ur] \ar@{-}[dr] & \;\; n \ge 4\\
& & & \stackrel{\alpha_{n}}{\bullet} } & \alpha_1 & \alpha_1 + \alpha_2 & 2\\
\hline
\xymatrix{\stackrel{\alpha_1}{\circ}\ar@{-}[r] <-.7ex> \ar@/^2ex/ @{<->}[rrr]&  \stackrel{\alpha_2}{\bullet} \ar@{.}[r] <-.7ex> & \stackrel{\alpha_{n-2}}{\bullet}  \ar@{-}[r] <-.7ex> & \stackrel{\alpha_n}{\circ} } \; n \ge 3& \alpha_1 \mbox{ ou } \alpha_n & \alpha_1 + \alpha_2 \mbox{ ou } \alpha_{n-1} + \alpha_n & 2\\
\hline
\xymatrix{\stackrel{\alpha_1}{\bullet} \ar@{-}[r] <-.7ex> &  \stackrel{\alpha_2}{\circ} \ar@{-}[r] <-.7ex> & \stackrel{\alpha_{3}}{\bullet}  \ar@{.}[r] <-.7ex> & \stackrel{\alpha_{n-1}}{\bullet} \ar@{<=}[r] <-.7ex> & \stackrel{\alpha_{n}}{\bullet} } \; n \ge 3 & \alpha_2 & \alpha_1 + \alpha_2 & 2\\
\hline
\xymatrix{\stackrel{\alpha_1}{\bullet} \ar@{-}[r] <-.7ex> & \stackrel{\alpha_2}{\bullet} \ar@{=>}[r] <-.7ex> &\stackrel{\alpha_3}{\bullet} \ar@{-}[r] <-.7ex> & \stackrel{\alpha_4}{\circ}  } & \alpha_4 & \alpha_3 + \alpha_4 & 1\\
\hline

\end{array}
$$ 
\end{center}
\end{dem}

{\bf Remarque :} Plus g\'en\'eralement, si $d$ n'est pas de la forme $l(\nu) + |J_\nu|$, pour un poids sp\'ecial $\nu$ (on note $J_\nu$ l'ensemble $\{1 \le i \le r \tq (\nu + \rho , \tilda{\alpha_i}) < 0 \}$), alors $H^d(X,\call{L}_\lambda) = 0$ pour tout poids sp\'ecial $\lambda$ de $\pic X$.   
\subsection{Cas d'\'egalit\'e}\label{sec:ega}

En g\'en\'eral, le th\'eorème \ref{th:orbf} ne suffit pas pour am\'eliorer les in\'egalit\'es (\ref{eq:ineg}). N\'eanmoins, nous allons voir un exemple de vari\'et\'e sym\'etrique complète $X$ pour lequel ce th\'eorème \ref{th:orbf} permet de transformer en \'egalit\'es les in\'egalit\'es (\ref{eq:ineg}) et donc de d\'eterminer tous les groupes de cohomologie de tous les fibr\'es en droites sur $X$. 

Cet exemple sera la compactification magnifique de l'espace sym\'etrique $Sp_{2n}/Sp_4 \croi Sp_{2n-4}$, pour $n \ge 4$. Cet espace sym\'etrique est appel\'e $C_n^{H,4}$ dans \cite[tab. 8]{loos} et son diagramme de Satake est :
\begin{equation}\label{eq:cnh4}
\xymatrix{
\stackrel{\alpha_1}{\bullet} \ar@{-}[r] <-.7ex> &  \stackrel{\alpha_2}{\circ} \ar@{-}[r] <-.7ex> & \stackrel{\alpha_{3}}{\bullet}  \ar@{-}[r] <-.7ex> & \stackrel{\alpha_{4}}{\circ} \ar@{-}[r] <-.7ex> & \stackrel{\alpha_{5}}{\bullet} \ar@{.}[r] <-.7ex>  & \stackrel{\alpha_{n-1}}{\bullet} \ar@{<=}[r] <-.7ex> & \stackrel{\alpha_n}{\bullet} }\p
\end{equation}

Il se trouve en effet que pour cet exemple de vari\'et\'e sym\'etrique complète $X$ on a une d\'ecomposition cellulaire qui \og s'adapte bien \fg\ \`a l'\'etude de la suite spectrale de Grothendieck-Cousin. 

Fixons d'abord quelques notations :

Soient $J$ la matrice $\left(\begin{array}{cc}
0 & I_n \\
-I_n & 0
\end{array}\right)$ et $J_0$ la matrice diagonale :
$$\left( \begin{array}{cccc}
-I_2 & & &\\
& I_{n-2} & &\\
&& I_2&\\
&&&I_{n-2}
\end{array}\right)
$$
(où $I_p$ est la matrice identit\'e de taille $p \croi p$).

Le groupe $G$ sera dor\'enavant le groupe $Sp_{2n}(\kk)$ \ie le sous-groupe des matrices $g \in \GL_{2n}(\kk)$ telles que : ${}^tg J g =J$.

L'automorphisme $\theta : G \to G$, $g \donne J_0gJ_0$ est une involution dont le groupe des points fixes v\'erifie : $G^\theta \iso Sp_4 \croi Sp_{2n-4}$. Nous noterons $X$ la compactification magnifique de $G/G^\theta$. C'est une vari\'et\'e sym\'etrique complète de rang $2$. En particulier, $X$ n'est pas de rang minimal car $\rang G - \rang G^\theta = 0$ ; le r\'esultat que nous \'enonçons ci-dessous (\cf le th\'eorème \ref{th:chq}) n'est donc pas contenu dans \cite{toi2}.

On choisit $\Phi,\Phi_0,\Phi_1, ...$ comme dans la section \ref{sec:debut}.

On a :

$$\Delta_1=\{\alpha_2,\alpha_4\}\;\;,$$

$$ \theta(\alpha_2) = -\alpha_1 - \alpha_2 -\alpha_3 \;\;,$$ $$\theta(\alpha_4) = \left\{ \begin{array}{ll}
-2\alpha_3 - \alpha_4 & \si n =4\\
-\alpha_3-\alpha_4 -2\alpha_5 -...-2\alpha_{n-1}-\alpha_n & \si n \ge 5
\end{array}\right. $$

(la num\'erotation des racines est donn\'ee par le diagramme (\ref{eq:cnh4})).
On a donc : $\pic X = \Z \omega_2 \plus \Z \omega_4$. 

\begin{thm}\label{th:chq}
Soit $X$ la compactification magnifique de l'espace symétrique $Sp_{2n}/Sp_{4} \croi Sp_{2n-4}$, $n \ge 4$. Pour tout faisceau inversible $\call{L}_\lambda$ sur $X$, de poids $\lambda \in \Z\omega_2 + \Z\omega_4$, on a :

$$H^d(X,\call{L}_\lambda) = \Plus_{J \sub \{2,4\}} \Plus_{{\nu \in \lambda + R_J \cap \Omega_J \atop \nu + \rho \mbox{ \scriptsize r\'egulier }} \atop l(\nu) +|J| = d} L(\nu^+)$$

(les notations sont celles du th\'eorème \ref{thm:euler}).

Plus concrètement :

$$H^d(X,\call{L}_\lambda) = 0$$
si $d \not\in \{0,5,4n-12,4n-8,4n-4,8n-21,8n-16\}$ 

et si $d= 0,5,4n-12,4n-8,4n-4,8n-21 \mbox{ ou } 8n-16$ alors :
$$H^d(X,\call{L}_\lambda) = \Plus_{\mu \in E^d_\lambda} L(\mu^+)$$
pour les ensembles de poids $E^d_\lambda$ qui sont donnés dans le tableau suivant :

$$\begin{array}{|c|c|}
\hline
d & E^d_\lambda \\
\hline
0 &  \left\{x\omega_2 + y\omega_4 \tq  x,y \ge 0\right\} \cap \lambda + \Z_{\le 0} \tilda{\alpha_2} + \Z_{\le 0}\tilda{\alpha_4}\\
\hline
5 & \{x\omega_2 + y \omega_4 \tq  x \le -4 ,  x+y \ge -2\} \cap \lambda + \Z_{>0} \tilda{\alpha_2} + \Z_{\le 0} \tilda{\alpha_4}\\
\hline 
4n-12 & \{x\omega_2+ y \omega_4 \tq   y \le -2n+5, x + 2y \ge -2n +5\} \cap \lambda +\Z_{\le 0} \tilda{\alpha_2} + \Z_{> 0}\tilda{\alpha_4} \\
\hline
4n-8 & \{x\omega_2 +y \omega_4 \tq  -2y - 2n +5 \le x \le -y -2n+ 3\} \cap \lambda + \Z_{>0} \tilda{\alpha_2} + \Z_{\le 0} \tilda{\alpha_4} \atop \uni \{x \omega_2 + y \omega_4 \tq  -y-2 \le x \le -2y-2n+1\} \cap \lambda + \Z_{\le 0} \tilda{\alpha_2} + \Z_{>0} \tilda{\alpha_4}\\
\hline
4n-4 &\{x \omega_2 + y \omega_4 \tq  y \ge 0, x+2y \le -2n+1 \} \cap \lambda +\Z_{>0}\tilda{\alpha_2} + \Z_{\le 0} \tilda{\alpha_4}\\
\hline
8n-21 & \{x \omega_2 + y\omega_4 \tq x \ge 0 ,   x+y \le -2n+3\} \cap \lambda + \Z_{\le 0}\tilda{\alpha_2} + \Z_{>0}\tilda{\alpha_4}\\
\hline
8n -16 & \{x\omega_2 + y \omega_4 \tq x \le -4, y \le -2n+5\} \cap \lambda + \Z_{>0}\tilda{\alpha_2} + \Z_{>0}\tilda{\alpha_4}\\
\hline
\end{array}$$

\end{thm}

{\bf remarques :} i) la dimension de $X$ est $8n-16$ ;

ii) si on pose pour tout $d \ge 0$ et tout $J \sub \{2,4\}$ :
$$\Omega_J^d := \left\{ \nu \in \Omega_J \tq \nu + \rho \mbox{ r\'egulier } \et l(\nu) + |J| = d \right\}$$
 alors on a :
$$E^d_\lambda = \bigcup_{J \sub \{2,4\}} \Omega_J^d \cap \lambda + R_J \p$$
 
\begin{figure}[h]
\begin{center}
\input{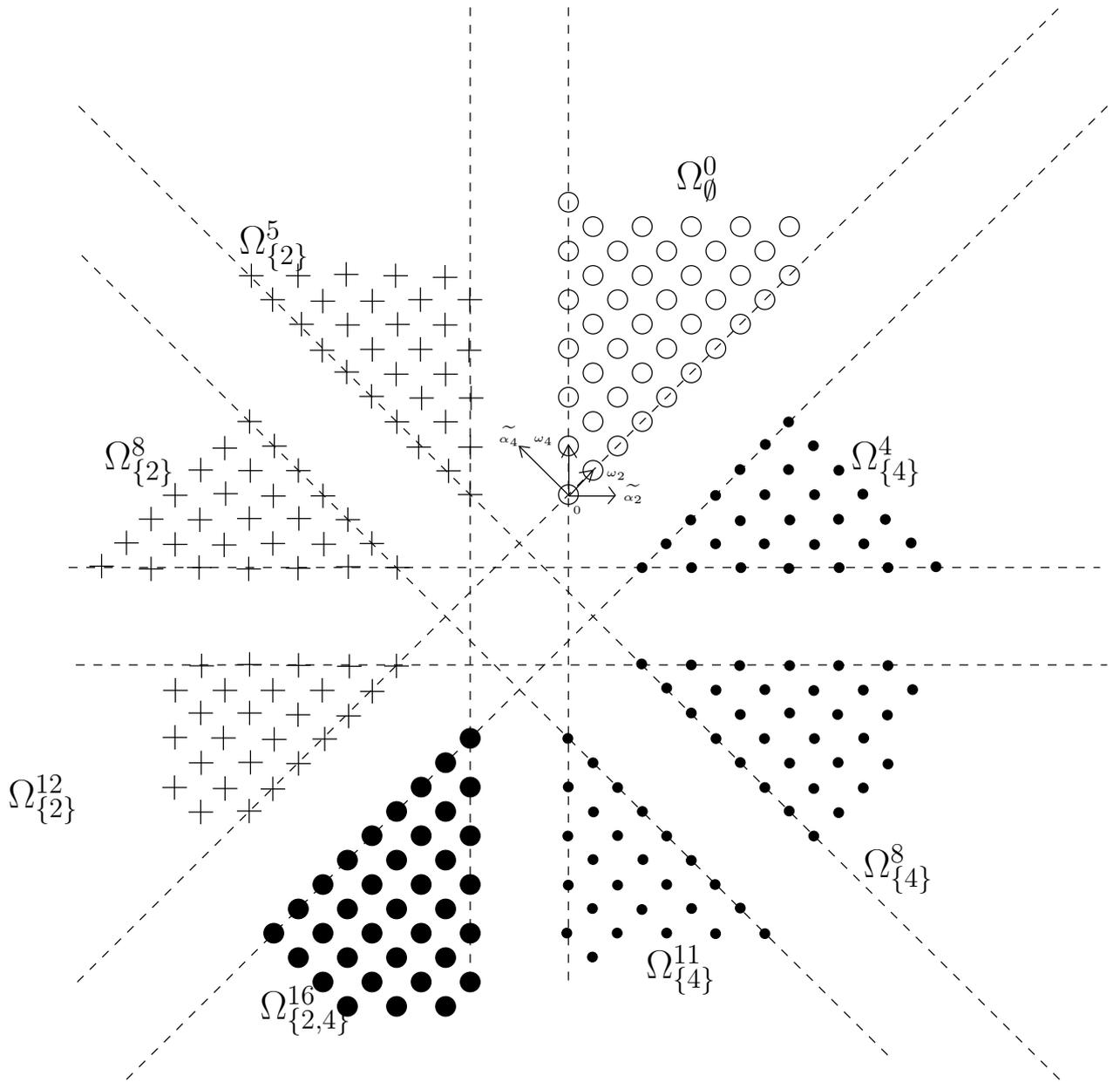}
    \caption{les ensembles $\Omega_J^d$ lorsque $n=4$}
\label{fig:ojd}
\end{center}
  \end{figure}

On a repr\'esent\'e les ensembles $\Omega_J^d$ dans le cas où $n=4$ sur la figure \ref{fig:ojd}.

\begin{dem}
Soit $\call{L}_\lambda$ un faisceau inversible sur $X$, de poids $\lambda \in \pic X$, tel que : $H^d(X,\call{L}_\lambda) \not= 0$. D'après le th\'eorème \ref{thm:ineg}, il existe alors un poids $\nu \in \pic X$ tel que : $\nu + \rho$ est r\'egulier et $l(\nu) + |J_\nu| = d$ (où $l(\nu) = |\{\alpha \in \Phi_1^+ \tq \cg \nu + \rho , \alpha^\ch \cd < 0 \}|$  et $J_\nu = \{i \in \{2,4 \} \tq (\nu + \rho , \tilda{\alpha_i}) < 0 \}$).

Or, on a :

$$\Phi_1^+ = \left\{ \sum_{i \le k < j} \alpha_k \tq i=1,2 \et j \ge 3 \mbox{ ou } i=3,4 \et j \ge 5 \right\}  $$
$$\dij \left\{ \sum_{i\le k <j} \alpha_k + 2 \sum_{j\le k < n} \alpha_k + \alpha_n \tq i=1,2,3,4 \et i < j \le n \right\} $$
$$\dij \left\{ 2 \sum_{i\le k \le n} \alpha_k + \alpha_n \tq i=1,2,3,4 \right\} \p$$

Pour un poids $\nu = x \omega_2 + y \omega_4 \in \pic X$, $x,y \in \Z$, lorsque $\alpha$ d\'ecrit l'ensemble $\Phi_1^+$, voici les valeurs prises par $\cg \nu + \rho , \alpha^\ch \cd $ :
 \begin{table}[htbp]\label{tab:nurho}
$$
\begin{array}{|c|c|}

\hline
\alpha \in \Phi_1^+ & {\cg \nu +\rho , \alpha^\ch \cd}\\
\hline

\ma \sum_{i \le k < j}\alpha_k & \left\{\begin{array}
{l}
x +j-i\si i=1,2 \et j = 3,4 \\
x+y + j-i\si i =1,2 \et 5 \le j \le n\\
  y +j-i\si i = 3,4 \et 5 \le j \le n
\end{array}\right. \\
\hline
\ma \sum_{i\le k < j} \alpha_k  + 2 \sum_{j \le k < n}\alpha_k +\alpha_n & \left\{\begin{array}{l}
2x + 2y +2n-1 \si i =1 ,j =2 \\
x + 2y +2n+2-i-j \si i =1,2 \et j =3, 4 \\
2y +2n -5 \si i =3 \et j =4\\
x+ y +2n+2-i-j\si i =1,2 \et 5 \le j \le n \\  

 y +2n+2-i-j \si i = 3,4 \et 5 \le j \le n
  
\end{array}\right.\\
\hline
\ma 2\sum_{i \le k < n}\alpha_k + \alpha_n & \left\{\begin{array}{l}x+y +n+1-i\si i =1,2 \\ y+n+1-i \si i =3,4
\end{array}\right. \\
\hline
\end{array}
$$
\caption{}
\end{table}

On déduit de ce tableau que le poids $\nu + \rho$ est régulier (\ie :$\cg \nu + \rho , \alpha^\ch \cd \not= 0 \cd$) si et seulement si les quatre conditions suivantes sont vérifiées :
$$\left\{\begin{array}{l}
x \not = -1,-2,-3\\
y\not= -1,-2, ..., -2n+6\\
x+y \not = -3,-4,...,-2n+4\\
x+2y \not = -2n+4  , -2n +3 , -2n +2 \p
\end{array}\right.$$

Puisque de plus :\begin{equation}\label{eq:J}
(\nu + \rho , \tilda{\alpha_2}) = c_n(x + 2) \et (\nu + \rho , \tilda{\alpha_4}) = c'_n(y+n-\frac{5}{2})
\end{equation}
(pour des constantes $c_n,c'_n >0$),
on obtient que $d=l(\nu) + |J_\nu|$ ne peut prendre que les valeurs $0, 5, 4n -12, 4n -8, 4n -4, 8n -21, 8n -16$ comme dans l'\'enonc\'e.
 
Pour terminer la d\'emonstration, il reste \`a v\'erifier que la \og composante finie \fg\ de la suite spectrale :
$$E_1^{p,q} = H^{p+q}_{X^+_p}(\call{L}_\lambda) \impliq H^{p+q}(X,\call{L}_\lambda)$$
d\'eg\'enère \ie : 
\begin{equation}\label{eq:ere1}
[E_r^{p,q} : L(\mu)] = [E_1^{p,q} : L(\mu)]
\end{equation}
pour tout poids dominant $\mu$ et pour tous $p,q \in \Z$ et tout $r \ge 1$.

On aura alors :
$$[H^d(X,\call{L}_\lambda) : L(\mu)] = \sum_{p,q \atop p+q = d}  [E_1^{p,q}:L(\mu)] \p$$

Fixons un poids dominant $\mu$.

Nous allons montrer (\ref{eq:ere1}) par r\'ecurrence sur $r$.

On suppose donc que $[E_r^{p,q} : L(\mu)] = [E_1^{p,q} : L(\mu)]$ et on consid\'ere les morphismes de $\goth g-$modules de la suite spectrale :
$$d_r^{p,q} : E_r^{p,q} \to E_r^{p+r,q-r+1} \p$$

Puisque $E_{r+1}^{p,q} = \ker  d_r^{p,q} / \im d_r^{p-r,q+r-1}$, il s'agit de d\'emontrer que pour tous $p, q \in \Z$, le $\goth g-$module $\im(d_r^{p,q})$ est de multiplicit\'e nulle selon $L(\mu)$. 

Raisonnons par l'absurde : supposons que :
\begin{equation}\label{eq:absu}
[\im d_r^{p,q} : L(\mu)] \not= 0 \p
\end{equation}

On a alors :
$$[E_r^{p,q} : L(\mu)] \not = 0 \et [E_r^{p+r,q-r+1} : L(\mu)] \not=0 \p$$

Comme $E_r^{p,q}$ est un sous-quotient de $E_1^{p,q}$, on a donc aussi :
$$[E_1^{p,q} : L(\mu)] = [H^{p+q}_{X^+_p}(\call{L}_\lambda): L(\mu)] \not=0 \p$$

En cons\'equence, d'une part, la cellule $X^+_p$ est de codimension $p+q$ et d'autre part, $X^+_p$ est centr\'ee en un point de l'orbite ferm\'ee $F$ de $X$ de la forme $x_p = w_\nu Q/Q$ pour un certain poids $\nu \in \pic X$ (d'après le th\'eorème \ref{th:orbf} et (\ref{eq:multcs})). D'où : $p+q = l(\nu) + |J_\nu|$ avec $\nu \in \pic X$ ; en particulier :\begin{equation}\label{eq:lJ}
p+q = 0 ,5 , 4n -12, 4n-8,4n-4 \mbox{ ou } 8n -16 \p
\end{equation}

De m\^eme, la cellule $X^+_{p+r}$ est centr\'ee en un point $x_{p+r}=w_{\nu'}Q/Q$ pour un certain poids $\nu ' \in \pic X$ et on aurait : \begin{equation}\label{eq:lJ2}
p+q+1 = 0 ,5 , 4n -12, 4n-8,4n-4 \mbox{ ou } 8n -16 \p
\end{equation} 
 
Il r\'esulte alors de (\ref{eq:lJ}) et (\ref{eq:lJ2}) que $1 = a-b$ pour certains $a,b \in \{0 ,5 , 4n -12, 4n-8,4n-4,8n -16 \}$ : cela est impossible si $n \ge 5$.

Il reste donc \`a traiter  le cas de la compactification magnifique de $Sp_8/Sp_4 \croi Sp_4$ \ie le cas où $n =4$.

Dans ce cas, il se peut que (\ref{eq:lJ}) et (\ref{eq:lJ2}) soient v\'erifi\'ees, on a donc besoin d'un argument suppl\'ementaire pour conclure. Les seuls cas possibles sont :
$$ (*) \left\{\begin{array}{l}
p+q = 4n -12 = 4\\
p+q+1 = 5
\end{array}\right. \et (**) \left\{\begin{array}{l} p+q = 8n -21 = 11 \\
p+q+1 = 4n-4 = 12 \p
\end{array}\right. $$

Ces deux cas se traitent de la m\^eme façon ; on va par exemple supposer $(*)$, c'est-\`a-dire que : $l(\nu) + |J_\nu| = p +q = 4$ et $l(\nu')+|J_{\nu'}| = 5$.

Comme $n=4$, le tableau 1 devient :

$$
\begin{array}{|c|c|}
\hline
\alpha \in \Phi_1^+ & {\cg \nu +\rho , \alpha^\ch \cd}\\
\hline
\ma \sum_{i \le k < j}\alpha_k & 
x +j-i\si i=1,2 \et j = 3,4 \\
\hline
\ma \sum_{i\le k < j} \alpha_k  + 2 \sum_{j \le k < 4}\alpha_k +\alpha_4 &\left\{\begin{array}{l}
2x + 2y +7 \si i =1 ,j =2 \\
x + 2y + 10 -i-j \si i =1,2 \et j =3, 4 \\
2y +3 \si i =3 \et j =4\\  
\end{array}\right.\\
\hline
\ma 2\sum_{i \le k < 4}\alpha_k + \alpha_n & \left\{\begin{array}{l}x+y +5-i\si i =1,2 \\ y+5-i \si i =3,4
\end{array}\right. \\
\hline
\end{array}
$$

On a : $l(\nu) + |J_\nu| = 4 \impliq |J_\nu| =1$. Mais alors, d'après (\ref{eq:J}), on a forc\'ement : $J_\nu = \{\tilda{\alpha_4}\}$. De m\^eme,  $l(\nu') + |J_{\nu'}| = 5 \impliq  J_{\nu'} = \{\tilda{\alpha_2}\} $.

En particulier, en posant $D:=D_{\tilda{\alpha_4}}$ :
$$X^+_p \sub D \et X^+_{p+r} \cap D \substr X^+_{p+r} \p$$

Soit ${\cal I}_D$ le faisceau d'id\'eaux d\'efinissant le diviseur $D$ (c'est un faisceau inversible). Le morphisme $\call{L}\tens {\cal I}_D \to \call{L}$ induit un morphisme entre les suites spectrales associ\'ees \`a $\call{L} \tens {\cal I}_{D}$ et \`a $\call{L}$. On a donc le diagramme commutatif de $\goth g-$modules suivant :
$$
\xymatrix{
E_r^{p,q}(\call{L}_\lambda) \ar^{d_r^{p,q}}[r] & E_r^{p+r,q-r+1}(\call{L}_\lambda) \\
E_r^{p,q}(\call{L}_\lambda \tens {\cal I}_D) \ar^{\alpha_r^{p,q}}[u] \ar^{d_{r,D}^{p,q}}[r] & E_r^{p+r,q-r+1}(\call{L}_\lambda \tens {\cal I}_D) \ar_{\alpha_r^{p+r,q-r+1}}[u]
}
$$

D'un autre c\^ot\'e, la suite exacte courte 
$$0 \to \call{L}_\lambda \tens {\cal I}_D \to \call{L}_\lambda \to \call{L}_\lambda\res{D} \to 0$$
induit une longue suite exacte :
$$... \to H^4_{X^+_p}(\call{L}_\lambda \tens {\cal I}_D) \to H^4_{X^+_p}(\call{L}_\lambda) \to H^4_{X^+_p}(\call{L}_\lambda \res{D}) \to ...$$

Si $$ X^+_p \sub D \;\;,$$ la cellule $X^+_p$ est de codimension $3$ dans $D$ et $ H^4_{X^+_p}(\call{L}_\lambda \res{D}) = 0$. Le morphisme :
$$H^4_{X^+_p}(\call{L}_\lambda \tens {\cal I}_D) = E_1^{p,q}(\call{L}\tens {\cal I}_D) \to E_ 1^{p,q}(\call{L}_\lambda) = H^4_{X^+_p}(\call{L}_\lambda)$$
est donc surjectif.

Cela n'entra\^inerait pas forc\'ement que le morphisme $\alpha_r^{p,q}$ du diagramme pr\'ec\'edent soit surjectif mais puisque, par hypothèse de r\'ecurrence, 
$$[E_1^{p,q}(\call{L}) : L(\mu)] = [E_r^{p,q}(\call{L}) : L(\mu)]$$
pour tout faisceau inversible $\call{L}$ sur $X$, on a n\'eanmoins :
$$[\im \alpha_r^{p,q}:L(\mu)] = [E_r^{p,q}:L(\mu)]$$
d'où :
$$[\im (d_r^{p,q}) : L(\mu)] = [\im (d_r^{p,q} \circ \alpha) : L(\mu)]\p$$

En cons\'equence :
$$[\im (d_r^{p,q}) : L(\mu)] = [\im (\alpha_r^{p+r,q-r+1} \circ d_{r,D}^{p,q}) : L(\mu)]$$
$$ \le [E_r^{p+r,q-r+1}(\call{L}\tens {\cal I}_D) : L(\mu)]$$
$$ \le [H^5_{X^+_{p+r}}(\call{L}_\lambda \tens {\cal I}_D) : L(\mu) ] \p$$

En raisonnant de m\^eme \`a partir du diagramme commutatif :$$
\xymatrix{
E_r^{p,q}(\call{L}_\lambda \tens {\cal I}_D^{n-1}) \ar[r] & E_r^{p+r,q-r+1}(\call{L}_\lambda \tens {\cal I}_D^{n-1}) \\
E_r^{p,q}(\call{L}_\lambda \tens {\cal I}_D^n) \ar[u] \ar[r] & E_r^{p+r,q-r+1}(\call{L}_\lambda \tens {\cal I}_D^n) \ar[u]
}
$$

on obtient que :

$$[\im (d_r^{p,q}) : L(\mu)] \le [H^5_{X^+_{p+r}}(\call{L}_\lambda \tens {\cal I}_D^n) : L(\mu) ] $$

pour tout $n \ge 1$.

Or, si $D \cap X^+_{p+r} \substr X^+_{p+r}$, on a :
$$\lim_{\longleftarrow \atop n} H^5_{X^+_{p+r}}(\call{L}_\lambda \tens {\cal I}_D^n) = 0 \p$$ 
On a donc : $[\im (d_r^{p,q}) : L(\mu)] = 0$ ;  d'où la contradiction (\cf (\ref{eq:absu})).

\end{dem}

\subsection{Cas de la vari\'et\'e des coniques complètes}

Nous allons d\'emontrer le r\'esultat suivant :
\begin{thm}\label{thm:ccc}
Soient $R_1:=2\Z_{> 0} \alpha_1 -2 \Z_{\ge 0}  \alpha_2$ et $R_2 := -2 \Z_{\ge 0} (\alpha_1 +  2\Z_{> 0} \alpha_2$. On note $s_1,s_2$ les r\'eflexions simples $s_{\alpha_1}$ et $s_{\alpha_2}$ et $\rho=\alpha_1+\alpha_2$.

Pour tout $\lambda  \in 2\Z\omega_1 + 2\Z\omega_2$, on a :

i) $$   H^0(\CCC,\call{L}_\lambda) =  \Plus_{\mu \in \lambda - 2\Z_{\ge 0} \alpha_1  - 2\Z_{\ge 0} \alpha_2 \atop \mu \; \mathrm dominant} L(\mu) \et H^5(\CCC,\call{L}_\lambda) = \Plus_{\mu \in w_0 \lambda - 2\Z_{\ge 4} \alpha_1  - 2\Z_{\ge 4} \alpha_2 \atop \mu \; \mathrm dominant} L(\mu) $$

ii)  $$H^1(\CCC,\call{L}_\lambda) = H^4(\CCC,\call{L}_\lambda) =0$$

iii) $$H^2(\CCC,\call{L}_\lambda) = \Plus_{\mu \; \mathrm{dominant} \atop \mu + \rho  \in s_1 (\lambda +\rho + R_1)  \setminus (s_2s_1)(\lambda +\rho +R_1) } L(\mu) \plus \Plus_{\mu \; \mathrm{dominant} \atop \mu + \rho \in s_2 (\lambda + \rho +R_2) \setminus (s_1s_2)(\lambda +\rho +R_2) }L(\mu) $$

$$H^3(\CCC,\call{L}_\lambda) =\Plus_{\mu \; \mathrm{dominant} \atop \mu \in (s_2s_1) (\lambda +\rho +R_1)  \setminus s_1(\lambda +\rho + R_1) } L(\mu) \plus \Plus_{\mu \; \mathrm{dominant} \atop \mu \in (s_1s_2) (\lambda +\rho + R_2) \setminus s_2(\lambda +\rho + R_2) }L(\mu) $$
\end{thm}
\begin{dem}

On va montrer iii). On traite le cas de $H^2(\CCC,\call{L}_\lambda)$ :

Nous verrons d'abord que si $\mu + \rho \in s_1(\lambda+\rho +R_1)\moins (s_2s_1)(\lambda+\rho +R_1)$ ou si $\mu \in s_2(\lambda+\rho +R_2)\moins (s_1s_2)(\lambda+\rho +R_2)$ alors $L(\mu)$ appara\^it dans la d\'ecomposition de $H^2(\CCC,\call{L}_\lambda)$ avec la multiplicit\'e $1$ puis que $\left[H^2(\CCC,\call{L}_\lambda):L(\mu)\right] = 0$ dans tous les autres cas.

Remarquons que si $\mu + \rho \in s_1(\lambda + \rho ++ R_1)$, il existe $m,n \in \Z$ tels que :
$$\mu = s_1(\lambda + \rho + 2m \alpha_1 + 2n \alpha_2) - \rho $$
$$ = \lambda + (2m-\cg \lambda,\alpha_1^\ch \cd -1)\alpha_1 + 2n \alpha_2 \p$$

Or, $\lambda \in \pic(\CCC) = 2 \Z \omega_1 + 2 \Z \omega_2$ donc $\cg \lambda,\alpha_1^\ch \cd \in 2\Z$. Donc $\mu \in \lambda + (2\Z+1) \alpha_1 + 2\Z \alpha_2$. 

De m\^eme, si $\mu +\rho \in s_2(\lambda+\rho +R_2)$ ou si $\mu +\rho \in (s_1s_2)(\lambda+\rho + R_2)$, on a $\mu \in \lambda+2\Z\alpha_1 + (2\Z+1)\alpha_2$.

En particulier, on a :
$$ (*) \;\; s_1(\lambda + \rho + R_1) \cap s_2(\lambda+\rho R_2) = s_1(\lambda + \rho +R_1) \cap (s_1s_2)(\lambda+\rho +R_2) = \vide \p$$

{\bf 1ère \'etape :} Si $\mu + \rho \in s_1(\lambda + \rho + R_1) \moins (s_2s_1)(\lambda +\rho +R_2)$, alors $\left[H^2(\CCC,\call{L}_\lambda) : L(\mu)\right] = 1$.

On se donne une d\'ecomposition de $\CCC$ en cellules de Bialynicki-Birula :
$$\CCC =\Dij_{x \in \CCC^T} \CCC^+_{x} \p$$

On a une suite spectrale convergente :$E_1^{p,q} \impliq H^{p+q}(\CCC,\call{L}_\lambda)$ pour une certaine num\'erotation $x_1,x_2,...$ des points fixes $x \in \CCC^T$ et où $E_1^{p,q}=H^{p+q}_{\CCC^+_{x_p}}(\call{L}_\lambda)$. Soit $p_0 \ge 1$ tel que $x_{p_0}=s_1\zzz$. Comme $X^+_{s_1} \sub D_{\alpha_1}$ et $X^+_{s_1}\not\sub D_{\alpha_2}$, on a d'après \cite[pro. 4.6 et th. 4.4]{toi} :
 $$\left[ E_1^{p_0,2-p_0}:L(\mu)\right] = \left[H^2_{X^+_{s_1}}(\call{L}_\lambda : L(\mu)\right] =1\p$$
En revanche, si $p\not= p_0$, on a pour tout $r \ge 0$ :
$$\left[ E_r^{p,2-p}:L(\mu) \right] \le \left[ E_1^{p,2-p}:L(\mu) \right] = 0 \p$$

En effet, sinon, d'une part on aurait d'après le th\'eorème \ref{th:orbf} $x_p$ de la forme $x_p = w\zzz$ pour un certain $w \in W$ tel que $\codim_{\CCC} \CCC^+_w =2$ et donc $x_p=s_2\zzz$ (car $x_p\not=x_{p_0}=s_1\zzz$). Mais d'autre part, on aurait aussi $\mu + \rho \in s_2(\lambda + \rho +R_2)$ ce qui est impossible d'après $(*)$.

On en d\'eduit que $[H^2(\CCC,\call{L}_\lambda):L(\mu)] = \left[E_\infi^{p_0,2-p_0} : L(\mu)\right]$.

Or, si $r \ge 1$, 
$$E_{r+1}^{p_0,2-p_0}= \frac{\ker (E_r^{p_0,2-p_0} \to E_r^{p_0+r,3-p_0-r})}{\im (E_r^{p_0-r,1-p_0+r} \to E_r^{p_0,2-p_0})} \p$$

Mais d'une part :
$$[E_r^{p_0-r,1-p_0-r}:L(\mu)] \le [E_1^{p_0-r,1-p_0-r}:L(\mu)] = [H^1_{\CCC^+_{x_{p_0-r}}}(\call{L}_\lambda) : L(\mu)] = 0$$
car si la cellule $\CCC^+_{x_{p_0-r}}$ est de codimension $1$, $x_{p_0-r}\not\in F$ ; d'autre part, on a aussi :
$$[E_r^{p_0+r,3-p_0+r}:L(\mu)] \le [E_1^{p_0+r,3-p_0+r}:L(\mu)]= [H^3_{\CCC^+_{x_{p_0+r}}}(\call{L}_\lambda) : L(\mu)] = 0$$
car si $x_{p_0+r}$ est de la forme $x_{p_0+r} = w\zzz$ pour un certain $w \in W$, et si la cellule correspondante $\CCC^+_w$ est de codimension $3$, on a $w = s_1s_2 \mbox{ ou } s_2s_1$ et $\mu +\rho \not\in w(\lambda+\rho +R_w)$ (\cf (*)).

On en d\'eduit que : $$\left[E_{r+1}^{p_0,2-p_0}:L(\mu)\right] = \left[E_r^{p_0,2-p_0} : L(\mu)\right]$$
pour tout $r \ge 1$ et donc que :
 $$\left[H^2(\CCC,\call{L}_\lambda):L(\mu)\right] = \left[E_\infi^{p_0,2-p_0} : L(\mu)\right]$$
$$= \left[E_1^{p_0,2-p_0}:L(\mu)\right] = 1 \p$$  

\vskip .5cm

De la m\^eme façon, on montre que si $\mu +\rho \in s_2(\lambda+\rho+R_2) \moins (s_1s_2)(\lambda+\rho +R_2)$, $[H^2(\CCC,\call{L}_\lambda):L(\mu)] = 1$.
\vskip .5cm
{\bf 2ème \'etape} Si par exemple\footnote{Cela peut arriver, {\it e.g.} : $\mu := \omega_2$, $\lambda := -14\omega_1+8\omega_2$} $\mu \in (s_2s_1)(\lambda + \rho +R_1) \cap s_1(\lambda+\rho+R_1)$, alors $[H^2(\CCC,\call{L}_\lambda):L(\mu)] = 0$.

Pour le d\'emontrer, on choisit la d\'ecomposition cellulaire de Bialynicki-Birula de $\CCC$ donn\'ee par le sous-groupe \`a un paramètre :
$$\nu : \kk^* \to T ,\; z \donne \left(\begin{array}{ccc}
z^3 & 0 & 0\\
0&z\inv & 0\\
0 & 0 & z^{-2}
\end{array}\right) \p$$

Pour ce sous-groupe \`a un paramètre, on a : $\CCC= \Dij_{p=0}^{11} \CCC^+_p$ avec en particulier, $\CCC^+_5= \CCC^+_{s_1}$ et $\CCC^+_6=\CCC^+_{s_2s_1}$.

On sait d\'ej\`a que :
$$\left[H^2(\CCC,\call{L}_\lambda):L(\mu)\right] \le \sum_{p=0}^{11} \left[E_r^{p,2-p} : L(\mu)\right]$$
pour tout $r \ge 1$. Or, comme les seuls points fixes $x_p \in \CCC^T$ tels que $x_p \in F$ et $\CCC^+_p$ soit de codimension $2$ dans $\CCC$ sont $s_1\zzz$ et $s_2\zzz$ et comme $\mu +\rho \not\in s_2(\lambda + \rho + R_2)$ (\cf (*)), on a forc\'ement :
$$\left[E_r^{p,2-p}:L(\mu)\right] \not=0 \impliq \left[H^2_{\CCC^+_p}(\call{L}_\lambda):L(\mu)\right] \not=0$$
$$\impliq x_p=s_1\zzz \impliq p=5 \p$$

Ainsi, pour tout $r \ge 1$, $\left[H^2(\CCC,\call{L}_\lambda):L(\mu)\right] \le \left[E_r^{5,-3}:L(\mu)\right]$.

On va montrer que $\left[E_2^{5,-3}:L(\mu)\right] = 0$.

Puisque $$E_2^{5,-3} =\frac{\ker (E_1^{5,-3} \stackrel{d_1^{5,-3})}{\longrightarrow} E_1^{6,-3}}{\im (E_1^{4,-3} \stackrel{d_1^{4,-3}}{\longrightarrow} E_1^{5,-3})} \;\;,$$ il suffit de v\'erifier que $\left[\ker d_1^{5,-3} : L(\mu)\right] = 0$.

Or, $d_1^{5,-3}$ est le morphisme 
$$d_1^{5,-3} : H^2_{\CCC^+_{s_1}}(\call{L}_\lambda) \to H^3_{\CCC^+_{s_2s_1}}(\call{L}_\lambda) $$ 
induit par la d\'ecomposition $\CCC^+_{s_1}=\CCC^+_{s_2s_1} \dij \CCC^+_{s_1} \moins \CCC^+_{s_2s_1}$ où $\CCC^+_{s_2s_1}$ est ferm\'e dans $\CCC^+_{s_1}$.

Pour tout $\goth g$ module $M$ et tout caractère central $\chi : Z(\goth g) \to \kk$ ($Z(\goth g)$ est le centre de l'algèbre enveloppante de $\goth g$) rappelons que :
$$M_\chi : =\left\{m \in M \tq \exist n >0 ,\, (\ker \chi)^n m = 0\right\}$$
est l'espace-propre g\'en\'eralis\'e associ\'e \`a $\chi$.

On note encore $\chi_\mu$ le caractère central par lequel $Z(\goth g)$ agit sur $L(\mu)$ et :
$$M[\mu] : =\left\{m \in M_{\chi_\mu} \tq \qq t \in T,\, t.m = \mu(t) m\right\}$$ le $T-$espace-propre de $M_{\chi_\mu}$ de poids $\mu$. 

On va utiliser que $[M:L(\mu)] = \dim_\kk M[\mu]$ pour tout $\goth g-{\chap{B}}-$module.

On sait que $H^2_{\CCC^+_{s_1}}(\call{L}_\lambda)[\mu]$ est de dimension $1$ autrement dit : 
$$H^2_{\CCC^+_{s_1}}(\call{L}_\lambda)[\mu] = \kk f_\mu $$
pour un certain \'el\'ement $f_\mu$ non nul de $H^2_{\CCC^+_{s_1}}(\call{L}_\lambda)$.

On a donc $\ker(d_1^{5,-3})[\mu] = \kk f_\mu \mbox{ ou } 0$.

Or, on peut d\'eterminer $f_\mu$ et montrer que $d^{5,-3}_1(f_\mu) \not=0$ (\cf les annexes).

On a donc $\ker(d_1^{5,-3})[\mu] = 0$ d'où : 
$ \left[\ker d_1^{5,-3} : L(\mu)\right] = 0$.

\vskip .5cm De m\^eme, si $\mu + \rho \in s_2(\lambda + \rho + R_2) \cap (s_1s_2)(\lambda + \rho + R_2)$, on peut montrer que  $[H^2(\CCC,\call{L}_\lambda):L(\mu)] = 0$. 

Enfin, si $\mu + \rho \not\in s_1(\lambda+ \rho + R_1) \cup s_2 (\lambda + \rho + R_2)$, on a $$[H^2(\CCC,\call{L}_\lambda) : L(\mu)] \le [H^2_{{\CCC}^+_{s_1}}(\call{L}_\lambda):L(\mu)] + [H^2_{{\CCC}^+_{s_2}}(\call{L}_\lambda):L(\mu)] = 0 \p$$

\vskip .5cm

Le cas de $H^3(\CCC,\call{L}_\lambda)$ se d\'emontre avec les m\^emes arguments.

Voici une cons\'equence du th\'eorème \ref{thm:ccc} :

\begin{cor}
Soient $\lambda_1,\lambda_2 \in \Z$ et $\lambda:=2\lambda_1 \omega_1 +2\lambda_2\omega_2 \in \pic (\CCC)$. On a :
$$H^0(\CCC,\call{L}_\lambda) \not=0 \equi \left\{\begin{array}{l}
\lambda_1+2\lambda_2 \ge 0\\
2\lambda_1 + \lambda_2 \ge 0
\end{array}\right.$$
$$H^2(\CCC,\call{L}_\lambda) \not=0 \equi \left\{\begin{array}{l}
\lambda_1 +3 \le 0 \\
\lambda_1+\lambda_2 +1 \ge 0
\end{array}\right. \mbox{ ou }\left\{\begin{array}{l}
\lambda_2 +3 \le 0 \\
\lambda_1+\lambda_2 +1 \ge 0
\end{array}\right. 
$$
$$ H^3(\CCC,\call{L}_\lambda) \not=0 \equi \left\{\begin{array}{l}
\lambda_1 -1 \ge 0 \\
\lambda_1+\lambda_2 +3 \le 0
\end{array}\right. \mbox{ ou }\left\{\begin{array}{l}
\lambda_2 -1 \ge 0 \\
\lambda_1+\lambda_2 +3 \le 0
\end{array}\right. $$
$$H^5(\CCC,\call{L}_\lambda) \not=0 \equi \left\{\begin{array}{l}
\lambda_1+2\lambda_2 +6 \le 0\\
2\lambda_1 + \lambda_2 + 6 \le 0
\end{array}\right.\p$$
\end{cor}

(\cf la figure \ref{fig:piccc})

Par exemple, si $\lambda = 12\omega_1 -6\omega_2$, $H^0(\CCC,\call{L}_\lambda)$ et $H^0(\CCC,\call{L}_\lambda)$ sont non nuls. Mais, pour tout $\lambda \in \pic \CCC$, on a $H^2(\CCC,\call{L}_\lambda)$ ou $H^3(\CCC,\call{L}_\lambda) = 0$.

\addcontentsline{toc}{section}{Annexes}

\section*{Annexes}

{\bf Calcul  de $[\ker (H^2_{\CCC^+_{s_1}}(\call{L}_\lambda) : L(\mu)]$ : }

Notons $S_3$ l'espace des matrices sym\'etriques $3 \croi 3$ \`a coefficients dans $\kk$ et $\pi : S_3 \moins\{0\} \croi S_3 \moins\{0\} \to \PP(S_3) \croi \PP(S_3)$ la surjection standard.

Si $\lambda = 2\lambda_1\omega_1 + 2\lambda_2 \omega_2 \in \pic \CCC$, $(\lambda_1,\lambda_2 \in \Z)$, alors :
$$\call{L}_\lambda = {\cal O}_{\PP(S_3)}(\lambda_1) \cten {\cal O}_{\PP(S_3)}(\lambda_2) \res{\CCC} \p$$
Donc, pour tout ouvert $V$ de $\CCC$, $\call{L}_\lambda(V)$ est l'espace des fonctions  $f$ r\'egulières sur $\pi\inv V$ telles que :
$$\qq s, s' \in \kk,\, \qq Q,Q' \in S_3,\, f(sQ,s'Q') = s^{\lambda_1}s'^{\lambda_2}f(Q,Q') \p$$

Nous allons introduire certaines fonctions rationnelles $X_i,Y_i$ sur $\CCC$. Pour cela rappelons la description de la cellule ouverte de $\CCC$ associ\'ee au sous-groupe \`a un paramètre :
$$\nu : \kk^* \to T ,\; z \donne \left(\begin{array}{ccc}
z^3 & 0 & 0\\
0&z\inv & 0\\
0 & 0 & z^{-2}
\end{array}\right) \;;$$

Pour tout $z \in \kk^*$ et pour tout $([Q],[Q']) \in \CCC$, on a :
$$\nu(z).([Q],[Q']) = \left({\left[\begin{array}{ccc}
z^{-6} Q_{1,1} & z^{-2}Q_{1,2} & z\inv Q_{1,3}\\
 z^{-2}Q_{1,2} & z^2 Q_{2,2} & z^3 Q_{2,3}\\
z\inv Q_{1,3} & z^3 Q_{2,3} & z^4Q_{3,3}
\end{array} \right],\left[\begin{array}{ccc}
z^{6} Q'_{1,1} & z^{2}Q'_{1,2} & z Q'_{1,3}\\
 z^{2}Q'_{1,2} & z^{-2} Q'_{2,2} & z^{-3} Q'_{2,3}\\
z Q'_{1,3} & z^{-3} Q'_{2,3} & z^{-4}Q'_{3,3}
\end{array} \right] }\right) \p$$

Puisque l'on a :
$$\zzz=\left({\left[\begin{array}{ccc}
1 & 0 & 0\\
0 & 0 &0\\
0&0&0
\end{array}\right],\left[\begin{array}{ccc}
0 & 0 & 0\\
0 & 0 &0\\
0&0&1
\end{array}\right]}\right)$$

la cellule ouverte $\CCC^+_0$ est d\'efinie par :
$$\CCC^+_0 = \left\{([Q],[Q']) \in \CCC \tq Q_{1,1}\not=0, Q'_{3,3}\not=0 \right\}\p$$

Cette cellule ouverte est isomorphe \`a l'espace affine $\A^5$ :
$$\A^5 \iso \CCC_0$$
$$ (x_i)_{1 \le i \le 5} \donne \left(\begin{array}{ccc}
1 & x_1 & x_3 \\
0 & 1 & x_2\\
0 & 0 & 1
\end{array}\right) \cdot \left({\left[\begin{array}{ccc}
1 & 0 & 0\\
0 & x_4 & 0\\
0 & 0 & x_4x_5
\end{array}\right],\left[ \begin{array}{ccc}
x_4x_5 & 0 & 0\\
0 & x_5 & 0\\
0 & 0 & 1
\end{array}\right]}\right)\p$$

Pour chaque $1 \le i \le 5$, on d\'efinit alors $X_i$ (respectivement $Y_i$) comme une fonction r\'egulière sur l'ouvert $s_1\CCC^+_0$ (respectivement $s_2s_1\CCC^+_0$) par :
$$X_i : s_1\cdot\left({\left(\begin{array}{ccc}
1 & x_1 & x_3 \\
0 & 1 & x_2\\
0 & 0 & 1
\end{array}\right) \cdot \left({\left[\begin{array}{ccc}
1 & 0 & 0\\
0 & x_4 & 0\\
0 & 0 & x_4x_5
\end{array}\right],\left[ \begin{array}{ccc}
x_4x_5 & 0 & 0\\
0 & x_5 & 0\\
0 & 0 & 1
\end{array}\right]}\right)}\right) \donne x_i $$
(respectivement : 
$$Y_i : (s_2s_1)\cdot\left({\left(\begin{array}{ccc}
1 & y_1 & y_3 \\
0 & 1 & y_2\\
0 & 0 & 1
\end{array}\right) \cdot \left({\left[\begin{array}{ccc}
1 & 0 & 0\\
0 & y_4 & 0\\
0 & 0 & y_4y_5
\end{array}\right],\left[ \begin{array}{ccc}
y_4y_5 & 0 & 0\\
0 & y_5 & 0\\
0 & 0 & 1
\end{array}\right]}\right)}\right)  \donne y_i \;\;).$$

Or, comme on a :
$$s_1\zzz=\left({\left[\begin{array}{ccc}
0& 0 & 0\\
0 & 1 &0\\
0&0&0
\end{array}\right],\left[\begin{array}{ccc}
0 & 0 & 0\\
0 & 0 &0\\
0&0&1
\end{array}\right]}\right)$$
$$s_2s_1\zzz=\left({\left[\begin{array}{ccc}
0 & 0 & 0\\
0 & 0 &0\\
0&0&1
\end{array}\right],\left[\begin{array}{ccc}
0 & 0 & 0\\
0 & 1 &0\\
0&0&0
\end{array}\right]}\right)\;,$$

les cellules $\CCC^+_{s_1}$ et $\CCC^+_{s_2s_1}$ sont donn\'ees par :
 $$\CCC^+_{s_1}= \left\{ ([Q],[Q']) \in \CCC \tq {Q_{1,1}=Q_{1,2}=Q_{1,3} =0, Q_{2,2}\not=0,\atop Q'_{3,3}\not=0} \right\}$$
$$\CCC^+_{s_2s_1}= \left\{([Q],[Q']) \in \CCC \tq {Q_{1,1}=Q_{1,2}=Q_{1,3} = Q_{2,2}= Q_{2,3} =0, Q_{3,3}\not=0 , \atop  Q'_{2,3} =Q'_{3,3}=0, Q'_{2,2}\not=0 }\right\}\p$$

On peut aussi v\'erifier que :
$$\CCC^+_{s_1} = \left\{ ([Q],[Q']) \in s_1 \CCC^+_0 \tq X_1=X_4 = 0 \right\}$$
$$\CCC^+_{s_2s_1}=\left\{ ([Q],[Q']) \in s_2s_1 \CCC^+_0 \tq Y_1=Y_3=Y_4 = 0 \right\} \p$$
En notant, pour $1\le i\le 5$, $x_i :=X_i \circ \pi$, fonction r\'egulière sur $\pi\inv\CCC^+_0$, et, pour $1 \le i,j \le 3$, $q_{i,j}$ (respectivement $q'_{i,j}$ la fonction r\'egulière sur $S_3 \croi S_3$ :
$$q_{i,j} : (Q,Q') \donne Q_{i,j} \mbox{ (respectivement }q_{i,j} : (Q,Q') \donne Q'_{i,j} \,)\;,$$
on a :
$$H^2_{\CCC^+_{s_1}}(\call{L}_\lambda) = \frac{\call{L}_\lambda(s_1\CCC^+_0 \cap (X_1 \not=0) \cap (X_4 \not=0))}{\call{L}_\lambda(s_1\CCC^+_0 \cap (X_1 \not=0)) + \call{L}_\lambda(s_1\CCC^+_0 \cap (X_4 \not=0))}$$
$$\iso \Plus_{n_2,n_3,n_5 \ge 0\atop n_1,n_4 >0}\kk \frac{x_2^{n_2}x_3^{n_3}x_5^{n_5}}{x_1^{n_1}x_4^{n_4}} q_{2,2}^{\lambda_1}{q'_{3,3}}^{\lambda_2}\p$$

Or, on a une r\'eunion d\'ecroissante :
$$H^ 2_{\CCC^+_{s_1}}(\call{L}_\lambda) = \Uni_{n \ge O} H^ 2_{\CCC^+_{s_1}}(\call{L}_\lambda \tens {\cal O}_{\CCC} (-nD_2) ) \p$$

Soient $M,N \ge 0$ deux entiers tels que :
\begin{eqnarray}\label{eq:MN}
\mu + \rho = s_1 (\lambda +\rho + 2(M+1) \alpha_1 -2N \alpha_2) \p
\end{eqnarray}

On a alors :

$$H^ 2_{\CCC^+_{s_1}}(\call{L}_\lambda \tens {\cal O}_{\CCC}(-nD_2)) = H^ 2_{\CCC^+_{s_1}}(\call{L}_\lambda) = \kk f_\mu$$
si $n \le N$ et $H^ 2_{\CCC^+_{s_1}}(\call{L}_\lambda \tens {\cal O}_{\CCC}(-nD_2)) = 0$ si $n > N$.

On en d\'eduit que :
$$f_\mu \in H^ 2_{\CCC^+_{s_1}}(\call{L}_\lambda \tens {\cal O}_\CCC(-ND_2)) \moins H^ 2_{\CCC^+_{s_1}}(\call{L}_\lambda \tens {\cal O}_\CCC(-(N+1)D_2)) \p$$
Comme $D_2 \cap s_1\CCC^+_0 = s_1 \CCC^+_0 \cap (X_5=0)$, cela signifie que :
$$f_\mu \in \Plus_{{n_2,n_3\ge 0 \atop n_5 \ge N}\atop n_1,n_4 >0}\kk \frac{x_2^{n_2}x_3^{n_3}x_5^{n_5}}{x_1^{n_1}x_4^{n_4}} q_{2,2}^{\lambda_1}{q'}_{3,3}^{\lambda_2} \moins \Plus_{{n_2,n_3\ge 0 \atop n_5 \ge N+1}\atop n_1,n_4 >0}\kk \frac{x_2^{n_2}x_3^{n_3}x_5^{n_5}}{x_1^{n_1}x_4^{n_4}} q_{2,2}^{\lambda_1}{q'}_{3,3}^{\lambda_2} \p$$

  Or, les fonctions $x_i$, $1 \le i \le 5$, $q_{2,2},{q'}_{3,3}$ sont des $T-$vecteurs propres de poids respectifs :
$$s_1 (-\alpha_1), s_1(-\alpha_2), s_1(-\alpha_1-\alpha_2), s_1(-2\alpha_1), s_1(-2\alpha_2), s_1(-2\omega_1), s_1(-2\omega_2) \p$$
Il s'ensuit que, pour $n_1,n_4 >0$, $n_2,n_3 \ge 0$ et $n_5 \ge N$, le mon\^ome $ \frac{x_2^{n_2}x_3^{n_3}x_5^{n_5}}{x_1^{n_1}x_4^{n_4}}q_{2,2}^{\lambda_1}{q'}_{3,3}^{\lambda_2}$ est de $T-$poids $\mu$ si et seulement si :
$$s_1(\lambda +(n_1+2n_4-n_3) \alpha_1 - (2n_5 +n_2 +n_3)\alpha_2) = \mu$$
$$\equi  \lambda +(n_1+2n_4-n_3) \alpha_1 - (2n_5 +n_2 +n_3)\alpha_2 = s_1\mu$$
$$\equi \lambda +(n_1+2n_4-n_3) \alpha_1 - (2n_5 +n_2 +n_3)\alpha_2= \lambda+(2M+3) \alpha_1 -2N\alpha_2$$
$$\equi \left\{\begin{array}{l}
n_1 +2n_4-n_3= 2M+3 \\
2(n_5-N)+n_2+n_3 =0
\end{array}\right.$$
$$\equi \left\{\begin{array}{l}
n_1 = 2(M-n_4)+3 \\
n_5=N,n_2=n_3 =0
\end{array}\right. \p$$

Ainsi :
$$f_\mu \in \Plus_{n_4 = 1}^{M+1} \kk \frac{x_5^N}{x_1^{2(M-n_4)+3}x_4^{n_4}}\p$$
De m\^eme, en  utilisant la filtration croissante :
$$H^2_{\CCC^+_{s_1}}(\call{L}_\lambda) = \uni_{m\ge 0}\ker \left(H^2_{\CCC^+_{s_1}}(\call{L}_\lambda)  \to H^2_{\CCC^+_{s_1}}(\call{L}_\lambda\tens {\cal O}_\CCC(mD_1)) \right)$$ on trouve aussi que :
$$f_\mu \in \Plus_{{n_2,n_3,n_5 \ge 0 \atop n_1 >0}\atop 1 \le n_4 \le M+1} \frac{x_2^{n_2} x_3^{n_3} x_5^{n_5}}{x_1^{n_1} x_4^{n_4}} q_{2,2}^{\lambda_1} {q'}_{3,3}^{\lambda_2} \moins \Plus_{{n_2,n_3,n_5 \ge 0 \atop n_1 >0}\atop 1 \le n_4 \le M} \frac{x_2^{n_2} x_3^{n_3} x_5^{n_5}}{x_1^{n_1} x_4^{n_4}}q_{2,2}^{\lambda_1}{q'}_{3,3}^{\lambda_2} \p$$

Finalement, \`a multiplication par un scalaire non nul près, $f_\mu$ est de la forme :
\begin{equation}\label{eq:fmu}
f_\mu = \frac{x_5^{N}}{x_1 x_4^{M+1}}q_{2,2}^{\lambda_1}{q'}_{3,3}^{\lambda_2} \end{equation}$$
+ \mbox{ une combinaison lin\'eaire des } \frac{x_5^{N}}{x_1^{2(M-n_4)+3} x_4^{n_4}}q_{2,2}^{\lambda_1}{q'}_{3,3}^{\lambda_2}$$
où $1 \le n_4 \le M$.

On peut montrer de la m\^eme façon que si $M',N'$ sont les entiers $\ge 0$ tels que :

\begin{eqnarray}\label{eq:M'N'}
\mu + \rho = (s_2s_1) (\lambda + \rho + 2(M'+1) \alpha_1 - 2N'\alpha_2) 
\end{eqnarray}
et si on pose $y_i:=Y_i \circ \pi$, fonctions r\'egulières sur $\pi\inv (\CCC^+_0)$, 
alors on a :
$$H^3_{\CCC^+_{s_2s_1}}(\call{L}_\lambda)[\mu] = \kk g_\mu$$
avec $g_\mu$ de la forme :

$$g_\mu = \frac{y_5^{N'}}{y_1 y_3 y_4^{M'+1}} q_{3,3}^{\lambda_1} {q'}_{2,2}^{\lambda_2} + \mbox{ une combinaison lin\'eaire des } \frac{y_2^{n_2} y_5^{n_5}}{y_1^{n_1} y_3^{n_3} y_4^{n_4}} q_{3,3}^{\lambda_1} {q'}_{2,2}^{\lambda_2}$$

où les $n_i$ sont des entiers v\'erifiant :

$$n_1 >0, n_2 \ge 0, n_3 > 0 , 0 < n_4 \le M'+1, n_5 \ge N' \et n_4 \le M' \mbox{ ou } n_5 > N'\p$$

Soient $\mu_1,\mu_2 \ge 0$ les entiers tels que : $\mu = \mu_1 \omega_1 + \mu_2 \omega_2$.

Remarquons que les relations (\ref{eq:MN}) et (\ref{eq:M'N'}) entra\^inent que \begin{eqnarray}\label{eq:MM'}
M' = M - \frac{\mu_2 +1}{2} \et N' = N +\frac{\mu_2 +1}{2}
\end{eqnarray}
(on sait d\'ej\`a que $\mu_2$ est impair).

Pour d\'eterminer $d_1^{5,-3}(f_\mu) \in H^3_{\CCC^+_{s_2s_1}}(\call{L}_\lambda)$, remarquons que la restriction :
$$\call{L}_\lambda \to \call{L}_\lambda\res{s_2s_1\CCC_0}$$
induit le diagramme commutatif suivant :
$$
\xymatrix{H^2_{\CCC^+_{s_1}}(\call{L}_\lambda) \ar[d] \ar[r]^{d_1^{5,-3}} & H^3_{\CCC^+_{s_2s_1}}(\call{L}_\lambda)\\
H^2_{\CCC^+_{s_1} \cap s_2s_1 \CCC_0}(\call{L}_\lambda) \ar[ru]_{d'} &.
}
$$

Il est plus facile de d\'ecrire le morphisme $d'$ que de d\'ecrire directement $d_1^{5,-3}$ ; on a en effet :
$$\CCC^+_{s_1}\cap s_2s_1 \CCC_0 = (s_2s_1 \CCC_0 \cap (Y_3 \not=0)) \cap (Y_2 = Y_4 = 0)$$
donc :
$$H^2_{\CCC^+_{s_1}\cap s_2s_1 \CCC_0}(\call{L}_\lambda) \iso H^1 (s_2s_1 \CCC_0 \cap (Y_3 \not=0) \moins (Y_2 = Y_4 = 0),\call{L}_\lambda)$$
$$\iso \frac{\call{L}_\lambda (s_2s_1 \CCC_0 \cap (Y_1Y_3Y_4 \not=0))}{\call{L}_\lambda (s_2s_1 \CCC_0 \cap (Y_1Y_3\not=0)) + \call{L}_\lambda (s_2s_1 \CCC_0 \cap (Y_3Y_4 \not=0))} $$
et le morphisme $d'$ est donn\'e par :
$$d' : H^2_{\CCC^+_{s_1}\cap s_2s_1 \CCC_0}(\call{L}_\lambda) \to H^3_{\CCC^+_{s_2s_1}}(\call{L}_\lambda)$$
$$ \qq n_2,n_5 \ge 0, n_1,n_4 >0, n_3 \in \Z,\; \frac{y_2^{n_2} y_5^{n_5}}{y_1^{n_1} y_3^{n_3} y_4^{n_4}} \donne \left\{\begin{array}{ll}
\frac{y_2^{n_2} y_5^{n_5}}{y_1^{n_1} y_3^{n_3} y_4^{n_4}} & \si n_1,n_4 < 0\\
0 & \sinon .
\end{array}\right.
$$

Maintenant, d\'ecrivons l'image de $f_\mu$ dans $H^2_{\CCC^+_{s_1}\cap s_2s_1 \CCC_0}(\call{L}_\lambda)$. Pour cela, on exprime les $x_i$ en fonction des $y_i$ dans le corps des fonctions rationnelles sur $S_3 \croi S_3$\footnote{
la notation \og $O(y_2)$ \fg signifie \og $y_2 \croi \mbox{ un polyn\^ome en les $y_i$ }$\fg
} :
$$
x_1 = \frac{y_1 y_3 + O(y_2)}{y_3^2+y_4y_5 + O(y_2)}$$
$$x_2= \frac{y_1y_5+ O(y_2)}{y_3^2+ y_4y_5 + y_1^2 y_5}$$
$$x_3= -\frac{y_3}{y_3^2+ y_4y_5 + y_1^2 y_5}$$
$$x_4=\frac{y_4(y_3^2+ y_4y_5 + y_1^2 y_5)}{(y_3^2+y_4y_5)^2 +O(y_2) }$$
$$x_5= \frac{y_5( y_3^2 + y_4y_5) + O(y_2)}{(y_3^2+ y_4y_5 + y_1^2 y_5)^2}
$$
(ces relations proviennent de l'\'egalit\'e \og matricielle \fg :
$$ 
s_1\cdot\left({\left(\begin{array}{ccc}
1 & x_1 & x_3 \\
0 & 1 & x_2\\
0 & 0 & 1
\end{array}\right) \cdot \left({\left[\begin{array}{ccc}
1 & 0 & 0\\
0 & x_4 & 0\\
0 & 0 & x_4x_5
\end{array}\right],\left[ \begin{array}{ccc}
x_4x_5 & 0 & 0\\
0 & x_5 & 0\\
0 & 0 & 1
\end{array}\right]}\right)}\right) $$
$$ = (s_2s_1)\cdot\left({\left(\begin{array}{ccc}
1 & y_1 & y_3 \\
0 & 1 & y_2\\
0 & 0 & 1
\end{array}\right) \cdot \left({\left[\begin{array}{ccc}
1 & 0 & 0\\
0 & y_4 & 0\\
0 & 0 & y_4y_5
\end{array}\right],\left[ \begin{array}{ccc}
y_4y_5 & 0 & 0\\
0 & y_5 & 0\\
0 & 0 & 1
\end{array}\right]}\right)}\right)
$$ sur $s_1 \CCC_0 \cap s_2s_1\CCC_0$).

D'un autre c\^ot\'e, comme on a les \'egalit\'es suivantes :
$$ \CCC^+_{s_1}\cap s_2s_1\CCC_0 = (s_2s_1\CCC_0 \cap s_1 \CCC_0) \cap (X_1=X_4 =0) $$
$$= (s_2s_1 \CCC_0\cap (Y_3 \not=0)) \cap (Y_1=Y_4 =0) $$
$$= (s_2s_1\CCC_0\cap s_1\CCC_0) \cap (X_4=X_1Y_1=0) \;,$$
on a les isomorphismes suivants (de $\goth g-$modules) :
$$H^2_{\CCC^+_{s_1}\cap s_2s_1 \CCC_0}(\call{L}_\lambda) \iso H^1((s_2s_1\CCC_0 \cap s_1 \CCC_0) \moins (X_1=X_4 =0),\call{L}_\lambda) $$
$$\iso H^1((s_2s_1 \CCC_0\cap (Y_3 \not=0)) \moins (Y_1=Y_4 =0),\call{L}_\lambda) $$
$$\iso H^1((s_2s_1\CCC_0\cap s_1\CCC_0) \moins (X_4=X_1Y_1=0),\call{L}_\lambda)
$$
$$\iso \frac{\call{L}_\lambda(s_2s_1\CCC_0 \cap s_1 \CCC_0\cap  (X_1X_4 \not=0))}{\call{L}_\lambda((s_2s_1\CCC_0 \cap s_1 \CCC_0) \cap (X_1\not=0)) + \call{L}_\lambda((s_2s_1\CCC_0 \cap s_1 \CCC_0) \cap (X_4 \not=0))} $$
$$\iso 
\frac{\call{L}_\lambda(s_2s_1\CCC_0 \cap (Y_1Y_3Y_4 \not=0))}{\call{L}_\lambda(s_2s_1\CCC_0 \cap  (Y_1Y_3\not=0)) + \call{L}_\lambda(s_2s_1\CCC_0 \cap (Y_3Y_4 \not=0))} $$
$$\iso
\frac{\call{L}_\lambda(s_2s_1\CCC_0 \cap s_1 \CCC_0\cap  (X_1X_4Y_1 \not=0))}{\call{L}_\lambda(s_2s_1\CCC_0 \cap s_1 \CCC_0 \cap (X_1Y_1\not=0)) + \call{L}_\lambda(s_2s_1\CCC_0 \cap s_1 \CCC_0 \cap (X_4 \not=0))} \p$$

En utilisant ces isomorphismes, en remplaçant dans la formule (\ref{eq:fmu}) de $f_\mu$ les $x_i$ par leur expression en les $y_i$, et en utilisant les relations (\ref{eq:MM'}), on trouve que l'image de $f_\mu$ dans $H^2_{\CCC^+_{s_1}\cap s_2s_1\CCC_0}(\call{L}_\lambda) $ est de la forme :
$$\frac{y_5^{N'}}{y_1y_3y_4^{M'+1}}q_{3,3}^{\lambda_1}{q'}_{2,2}^{\lambda_2}+\mbox{ une combinaison lin\'eaire de }\frac{y_2^{n_2}y_5^{n_5}}{y_1^{n_1}y_3^{n_3}y_4^{n_4}}q_{3,3}^{\lambda_1}{q'}_{2,2}^{\lambda_2}$$
où $n_5 > N'$ ou $n_2 >0$.

En cons\'equence : 
$$d^{5,-3}_1(f_\mu)\in H^3_{\CCC^+_{s_2s_1}}(\call{L}_\lambda) \iso \Plus_{n_2,n_5 \ge 0 \atop n_1,n_3, n_4 >0}\kk \frac{y_2^{n_2}y_5^{n_5}}{y_1^{n_1}y_3^{n_3}y_4^{n_4}} q_{3,3}^{\lambda_1}{q'}_{2,2}^{\lambda_2}$$
a un coefficient $1$ devant  $\frac{y_5^{N'}}{y_1y_3y_4^{M'+1}}q_{3,3}^{\lambda_1}{q'}_{2,2}^{\lambda_2}$.

Ainsi : $d^{5,-3}_1(f_\mu) \not=0$ (on a m\^eme $d^{5,-3}_1(f_\mu)=g_\mu$) et donc :
$$[\ker d_1^{5,-3} :L(\mu)] = 0\p$$
\end{dem}
\newpage

\newpage

\begin{figure}[h]
\begin{center}
\input{c2R.pstex_t}
    \caption{les ensembles $R_J$  pour la compactification magnifique de $\Sp_{2n}/\Sp_{4} \croi \Sp_{2n-4}$, $n \ge 5$}
    \label{fig:c2R}
\end{center}
  \end{figure}

\begin{figure}[h]\label{fig:piccc}
\begin{center}
\input{coniq1.pstex_t}
\caption{Le groupe de Picard de la vari\'et\'e des coniques complètes $\CCC$ (l\'egende : $\circ$ : $\lambda \in \pic \CCC$ tel que $H^0(\CCC,\call{L}_\lambda) \not=0$ ; $+$ : $\lambda \in \pic \CCC$ tel que $H^2(\CCC,\call{L}_\lambda) \not=0$ ; $\croi$ : $\lambda \in \pic \CCC$ tel que $H^3(\CCC,\call{L}_\lambda) \not=0$ ;  $\Box$  : $\lambda \in \pic \CCC$ tel que $H^5(\CCC,\call{L}_\lambda) \not=0$)}
\end{center}
\end{figure}

%%%BIBLIOGRAPHIE%%


\clearpage

\addcontentsline{toc}{section}{Bibliographie}

\begin{thebibliography}{BdCP99}


\bibitem[BGG]{BGG} I. N. {\sc Bernstein},  I. M. {\sc Gelfand},  S. I. {\sc Gelfand},
\textit{Differential operators on the base affine space and a study of ${\goth g}$-modules} in 
Lie groups and their representations (Proc. Summer School, Bolyai János Math. Soc., Budapest, 1971), pp. 21--64,
1975. 

\bibitem[Bialynicki-73]{BiBi} A. {\sc Bialynicki-Birula},
\textit{Some theorems on actions of algebraic groups},
 Annals of Mathematics, 98, pp. 480-497,
 1973.

\bibitem[Bialynicki-76]{BB76} A. {\sc Bialynicki-Birula},
\textit{Some properties of the decompositions of algebraic varieties determined by actions of a torus},
Bull. Acad. Pol. Sci., s\'er, sci. math. astron. phys. 24, pp.667-674, 1976.

\bibitem[Ch-Ma]{CHMA} R. {\sc Chirivi}, A. {\sc Maffei}, \textit{The ring of sections of a complete symmetric variety},  J. Algebra  261,  no. 2, pp. 310--326, 2003.

\bibitem[DeConcini-Pro\-ce\-si]{DCP} C. {\sc De Concini}, C. {\sc Procesi},
\textit{Complete symmetric varieties},
Invariant theory, Proc. 1st 1982 Sess. C.I.M.E., Montecatini/Italie, Lect. Notes Math. 996,
pp 1-44,
1983. 


%\bibitem[DeConcini-Springer]{DCS} C. {\sc De Concini}, T. A. {\sc Springer},
%\textit{Betti numbers of complete symmetric varieties}, Giornate di geometria, Roma 1984, in Geometry of today, PM 60, Birkh\"{a}user, 1985.

\bibitem[Dixmier]{Dix} J. {\sc Dixmier}, \textit{Algèbres enveloppantes}, Gauthiers-Villars, 1974.

\bibitem[Kempf]{Kempf} G. {\sc Kempf},
\textit{ The Grothendieck-Cousin complex of an induced representation},
Advances in Mathematics 29, pp 310-396,
1978.

\bibitem[Loos]{loos} O. {\sc Loos} Symmetric spaces II,  W. A. Benjamin, 1969.

\bibitem[Steinberg]{St} R. {\sc Steinberg} \textit{G\'en\'erateurs, relations et rev\^etements de groupes alg\'ebriques, pp. 113-127, colloq. th\'eorie des groupes alg\'ebriques de Bruxelles, Gauthiers-Villars, 1962}

\bibitem[T0]{these} A. {\sc Tchoudjem} \textit{Repr\'esentations d'algèbres de Lie dans des groupes de cohomologie \`a support}, Thèse de doctorat, Grenoble 2002 (disponible \`a http://www-fourier.ujf-grenoble.fr/THESE/ps/t125.ps.gz)

\bibitem[T1]{toi} A. {\sc Tchoudjem}
\textit{Cohomologie des fibrés en droites sur les compactifications des groupes réductifs},
Ann. Sci. École Norm. Sup. (4) 37, no. 3, 415--448,
2004.


\bibitem[T2]{toi2} A. {\sc Tchoudjem}
\textit{Cohomologie des fibrés en droites sur les vari\'et\'es magnifiques de rang minimal},
Bull. soc. math. de France, 135 (2), pp. 171-214, 2007.



%\bibitem[Wasserman]{Was} B. {\sc Wasserman},
%\textit{Wonderful varieties of rank two},
%Transformation groups, Vol. 1, no 4, pp. 375-403,
%1996.

 \end{thebibliography}
\end{document}